\documentclass[A4j,11pt]{article}
\usepackage{graphics}
\usepackage{amsmath}
\usepackage{amssymb}
\usepackage{latexsym}
\usepackage{amsfonts}
\begin{document}

\centerline{{\huge\bf A holonomy invariant}}
\centerline{{\huge\bf anisotropic surface energy}}
\centerline{{\huge\bf in a Riemannian manifold}}

\vspace{1truecm}

\centerline{\Large Naoyuki Koike}

\vspace{0.5truecm}

\begin{abstract}
In this paper, we investigate a holonomy invariant elliptic anisotropic surface energy for 
hypersurfaces in a complete Riemannian manifold, where ``holonomy invariant" means that the elliptic 
parametric Lagrangian (i.e., a Finsler metric) of the Riemannian manifold used to define the 
anisotropic surface energy is constant along each holonomy subbundle of the tangent bundle of 
the Riemannian manifold.  First we obtain the first variational formula for this anisotropic surface 
energy.  Next we shall introduce the notions of an anisotropic convex hypersurface, 
an anisotropic equifocal hypersurface and an anisotropic isoparametric hypersurface for this anisotropic 
surface energy.  
Also, we shall introduce the notion of an anisotropic tube for this anisotropic surface energy.  
We prove that anisotropic tubes over a one-point set in a symmetric space are anisotropic convex 
hypersurfaces and that anisotropic tubes over a certain kind of reflective submanifold 
in a symmetric space are anisotropic isoparametric and anisotropic equifocal hypersurfaces.  
\end{abstract}

%


\vspace{0.4truecm}

\noindent
{\Large\bf Introduction}

\vspace{0.3truecm}

H. Federer ([F1,2,3]) studied the elliptic parametric functional given by 
a parametric Lagrangian of general degree in the Euclidean space 
from the point of wide view of geometric measure theory, where we note that his study can be apply to 
the study of the elliptic parametric functional for (smooth) submanifolds of general codimension 
in the Euclidean space because the parametric Lagrangian is of general degree.  
B. White ([W]) studied the elliptic parametric functional for (smooth) submanifolds 
(of general codimension) in the Euclidean space in detail.  
In particular, U. Clarenz ([Cl]) studied the elliptic parametric functional for (smooth) hypersurfaces 
in the Euclidean space in more detail.  
On the other hand, Koiso and Palmer ([KP1-3, Palm]) studied a special elliptic parametric functional 
(which they called an anisotropic surface energy) for (smooth) hypersurfaces in the Euclidean space.  

In the case where the ambient space is a general complete oriented Riemannian manifold, 
Lira and Melo ([LM]) have recently introduced an elliptic parametric functional for (smooth) 
hypersurfaces as follows.  
Let $\widetilde M$ be an $(n+1)$-dimensional complete oriented Riemannian manifold and $M$ be 
an $n$-dimensional compact oriented manifold, where $M$ may have the boundary.  
Let $\widetilde F$ be a positive $C^{\infty}$-function over the tangent bundle 
$T\widetilde M\setminus\{0\}$ of 
$\widetilde M$ satisfying the homogeneity condition 
$$\widetilde F(tX)=t\widetilde F(X)\quad\,\,(X\in T\widetilde M\setminus\{0\},\,\,t>0).\leqno{({\rm H})}$$
Then $\widetilde F$ is called a {\it parametric Lagrangian} (of $\widetilde M$).  
Furthermore, if it satifies the elliptic condition 
$$(\widehat{\nabla}d\widetilde F)_X(Y,Y)>0\,\,\quad\,\,
(X\in T\widetilde M\setminus\{0\},\,\,Y(\not=0)\in T_X(T\widetilde M)\,\,{\rm s.t.}\,\,
\widehat g(X,Y)=0),
\leqno{({\rm E})}$$
then it is said to be {\it elliptic}, 
where $\widehat g$ is the Sasaki metric of $T\widetilde M$ and 
$\widehat{\nabla}$ is the Riemannian connection of $\widehat g$.  
Note that an elliptic parametric Lagrangian means a Finsler metric.  
Let $f$ be an immersion of $M$ into $\widetilde M$.  
They studied the following type of functional:
$${\cal F}(f):=\int_{x\in M}\widetilde F(\xi_x)dV,$$
where $\xi$ is the unit normal vector field of $f$ (compatible with the orientations of $M$ and 
$\widetilde M$) and $dV$ is the volume element of the induced metric on $M$ by $f$.  
They called ${\cal F}$ an {\it elliptic parametric functional}.  
In this paper, we shall call ${\cal F}$ an {\it anisotropic surface energy} following to 
the terminology of [KP1-3].  
In particular, if $\widetilde F$ is constant along each holonomy subbundle of $T\widetilde M$ 
(i.e., horizontally constant in the sense of [LM]), then 
we shall say that ${\cal F}$ and $\widetilde F$ are {\it holonomy invariant}.  
Note that the anisotropic surface energy treated by Koiso and Palmer is holonomy invariant.  
Lira and Melo ([LM]) stated that the horizontally constancy (=holonomy invariantness) of the elliptic 
anisotropic surface energy need to be imposed in order that 
the anisotropic mean curvature (which is a variational notion) of $f$ is given as 
the trace of the anisotropic analogue (which is a hypersurface theoretical notion) 
of the shape operator of $f$.  

\vspace{0.5truecm}

\noindent
{\bf Motivation 1.} 
As in the above statement by Lira and Melo, we consider that the elliptic anisotropic surface energy 
need to be assumed to be holonomy invariant in order to study the variational problem for the above 
elliptic anisotropic surface energy from the view of point of the hypersurface theory.  

\vspace{0.5truecm}

\noindent
{\bf Motivation 2.} 
In the case where the ambient space is a standard space, it is important to give 
examples of critical points (i.e., hypersurfaces with constant anisotropic mean curvature) of 
the holonomy invariant elliptic anisotropic surface energy over the class of all volume-preserving 
variations.  

\vspace{0.5truecm}

If the holonomy group of the ambient space $\widetilde M$ is big, then 
the holonomy invariantness of the elliptic anisotropic surface energy is a strong constraint condition.  
In particular, if $\widetilde M$ is rotationally symmetric (for example, a rank one symmetric space), 
then the only holonomy invariant elliptic anisotropic surface energy is a constant-multiple of 
the volume functional.  Hence we are interested in the case where the holonomy group of 
the ambient space is small.  
%

J. Ge and H. Ma [GM] studied the anisotropic surface energy treated by Koiso-Palmer.  
For this anisotropic surface energy, they introduced the notions of an anisotropic principal curvature 
and an anisotropic parallel translation of a hypersurface in the Euclidean space.  They proved that 
a hypersurface is with constant anisotropic principal curvatures if and only if 
the anisotropic parallel hypersurfaces of the hypersurface are of constant anisotropic 
mean curvature (i.e., the hypersurface is anisotropic isoparametric in the sense of [GM]).  
Furthermore, they obtained a Cartan identity for a hypersurface with 
constant anisotropic principal curvatures and classified complete 
hypersurfaces with constant anisotropic principal curvatures in terms of this identity.  

\vspace{0.5truecm}

\noindent
{\bf Motivation 3.} 
In a (general) complete Riemannian manifold, 
we should introduce the notions similar to the above notions introduced in [GM] and 
study the relations between the notions.  Also, we should give examples of the notions 
in a standard complete Riemannian manifold.  

\vspace{0.5truecm}

Under the above motivations, in this paper, we shall first give the first variational formula for 
a holonomy invariant elliptic anisotropic surface energy satisfying some convexity condition 
in a complete Riemannian manifold, which is a special one of the first variational formula for 
a (not necessarily holonomy invariant) elliptic anisotropic surface energy given by Lira and Melo ([LM]) 
but cannot be derived directly from their formula (see Section 1).  
Furthermore, by using this formula, we investigate a critical point of this anisotropic surface energy 
over the class of all variations (or over that of all volume-preserving variations) (see Section 1).  
Next, for this anisotropic surface energy, we introduce the notions of an anisotropic convex hypersurface, 
an anisotropic equifocal hypersurface and an anisotropic isoparametric hypersurface (see Section 2).  
Next, for this anisotropic surface energy, we introduce the notion of an anisotropic tube over a submanifold (Section 3) 
and prove that anisotropic tubes over an one-point set (which are called anisotropic geodesic spheres and 
are the notion corresponding to the Wulff shape in the Euclidean case) in a symmetric space of 
compact type or non-compact type are anisotropic convex hypersurfaces (see Section 4) and that 
anisotropic tubes over a certain kind of reflective submanifold in a symmetric space of compact type 
or non-compact type are anisotropic isoparametric and anisotropic equifocal hypersurfaces 
(see Section 5).  
Finally, we prove that the anisotropic equifocality is equivalent to the anisotropic isoparametricity 
in the case where the ambient space is a symmetric space of non-negative curvature (see Section 6).  

\section{Holonomy invariant anisotropic surface energy} 
In this section, we shall define a holonomy invariant anisotropic surface energy 
in a complete Riemannian manifold, which is treated in this paper, and obtain the first 
variational formula for this anisotropic surface energy.  

Let $\widetilde M$ be an $(n+1)$-dimensional complete (oriented) Riemannian manifold.  
Denote by $\langle\,\,,\,\,\rangle$ and $\widetilde{\nabla}$ the (Riemannian) metric and 
the Riemannian connection of $\widetilde M$, repsectively.  
Also, denote by $\tau_c$ the parallel translation along a curve $c$ in $\widetilde M$, 
$\Phi_p$ the holonomy group of $\widetilde M$ at $p(\in\widetilde M)$ and 
$P^{\rm hol}_v$ the holonomy subbundle of $T\widetilde M$ through $v(\in T\widetilde M)$.  
Also, denote by $S^n(1)_p$ the unit sphere centered at the origin in $T_p\widetilde M$.  
Fix $p_0\in\widetilde M$.  For simplicity, set $\Phi:=\Phi_{p_0}$ and $S^n(1):=S^n(1)_{p_0}$.  
It is clear that $\Phi$ acts on $S^n(1)$.  
Denote by $\pi_{\Phi}$ the orbit map of this action.  
For each  $p\in\widetilde M$, take a shortest geodesic $\gamma_p:[0,1]\to\widetilde M$ 
connecting $p_0$ to $p$ (i.e., $\gamma_p(0)=p_0,\,\gamma_p(1)=p$).  
Note that the choice of $\gamma_p$ is not unique for each $p$ belonging to the cut locus 
(which is denoted by $C$) of $p_0$.  In the sequel, we fix the choices of $\gamma_p$'s ($p\in C$).  
For simplicity, set $\tau_p:=\tau_{\gamma_p}$.  
Let $M$ be a $n$-dimensional compact (oriented) manifold, which may have boundary, 
and $f$ an immersion of $M$ into $\widetilde M$.  Denote by $\xi$ 
the unit normal vector field of $f$ (compatible with the orientations of $M$ and $\widetilde M$).  
Define a map $\nu:M\to S^n(1)$ by 
$$\nu(x)=\tau_{f(x)}^{-1}(\xi_x)\,\,\,\,(x\in M).$$
Set $\overline{\nu}:=\pi_{\Phi}\circ\nu$.  
Note that $\nu$ depends on the choices of $\gamma_p$'s ($p\in C$) but 
$\overline{\nu}$ is independent of their choices.  
Under fixed choices of $\gamma_p$'s ($p\in C$), we call $\nu$ the {\it Gauss map} of $f$.  
Denote by $C^{\infty}(S^n(1))_{\Phi}$ the ring of all $\Phi$-invariant $C^{\infty}$-functions 
over $S^n(1)$.  Take an elliptic parametric Lagrangian $\widetilde F$ of $\widetilde M$.  
Assume that $\widetilde F$ is holonomy invariant, that is, 
$\widetilde F\vert_{P^{\rm hol}_v}$ is constant for any $v\in T\widetilde M\setminus\{0\}$.  
Denote by $F_p$ the restriction of $\widetilde F$ to $S^n(1)_p$.  In particular, we denote 
$F_{p_0}$ by $F$ for simplicity.  
Then, since $\widetilde F$ is holonomy invariant, $F_p=F\circ\tau_p^{-1}$ holds for any 
$p\in\widetilde M$ and $F$ is $\Phi$-invariant.  
Thus $\widetilde F$ is determined by $F$, that is, any holonomy invariant elliptic parametric Lagrangian 
(of $\widetilde M$) is constructed from a $\Phi$-invariant $C^{\infty}$-function over $S^n(1)$.  
Furthermore, assume that $F$ satisfies the following {\it convexity condition}:
$$\nabla^S{\rm grad}\,F+F{\rm id}\,>\,0,
\leqno{({\rm C})}$$
where $\nabla^S$ is the Riemannian connection of $S^n(1)$, 
${\rm id}$ is the identity transformation of $TS^n(1)$.  
Note that the orbit space $S^n(1)/\Phi$ is an $(r-1)$-dimensional orbifold in the case where 
$\widetilde M$ is a symmetric space of rank $r$.  In particular, if $\widetilde M$ is 
a symmetric space of rank one, then $S^n(1)/\Phi$ is of dimension zero and hence $F$ must be 
constant and hence the above anisotropic surface energy ${\cal F}$ is equal to the (usual) volume functional 
up to a constant-multiple.  
Thus, if $\widetilde M$ is a symmetric space, then we are interesting 
in the case where the rank of the symmetric space is higher.  
Denote by ${\rm Imm}(M,\widetilde M)$ the set of all ($C^{\infty}$-)immersions of $M$ into 
$\widetilde M$.  
In this paper, we shall investigate the holonomy invariant elliptic anisotropic surface energy 
${\cal F}:{\rm Imm}(M,\widetilde M)\to{\Bbb R}$ by 
$${\cal F}(f):=\int_{x\in M}\widetilde F(\xi_x)dV,$$
where $dV$ is the volume element of the induced metric on $M$ by $f$.  
Note that $\widetilde F(\xi_x)=(F\circ\nu)(x)$ ($x\in M$) and hence $F\circ\nu$ is independent of 
the choices of $\gamma_p$'s ($p\in C$).  
Denote by $H$ the mean curvature of $f$ (with respect to $\xi$).  
Define a function $H_F$ over $M$ by 
$$(H_F)_x:=(F\circ\nu)(x)H_x-({\rm div}(f_{\ast}^{-1}(\tau_{f(\cdot)}
(({\rm grad}F)_{\nu(\cdot)}))))_x\,\,\,\,(x\in M),\leqno{(1.1)}$$
where ${\rm div}(\cdot)$ is the divergence of $(\cdot)$ with respect to the induced metric 
on $M$ by $f$.  Here we note that $\tau_{f(y)}(({\rm grad}F)_{\nu(y)})\in f_{\ast}(T_yM)$ 
for any $y\in M$ (see Figure 1).  We call $H_F$ the {\it anisotropic mean curvature} of $f$.  
If $H_F=0$, then we call $f:M\hookrightarrow\widetilde M$ an 
{\it anisotropic minimal hypersurface}.  
Also, if $H_F$ is constant, then we call 
$f:M\hookrightarrow\widetilde M$ a {\it hypersurface with constant anisotropic mean curvature}.  

\vspace{0.2truecm}

\centerline{
\unitlength 0.1in
\begin{picture}( 51.5800, 17.2500)(-14.3000,-20.9400)
%
\special{pn 8}%
\special{pa 1178 1126}%
\special{pa 738 1730}%
\special{pa 1616 1730}%
\special{pa 2052 1126}%
\special{pa 2052 1126}%
\special{pa 1178 1126}%
\special{fp}%
%
\special{pn 20}%
\special{sh 1}%
\special{ar 1350 1508 10 10 0  6.28318530717959E+0000}%
\special{sh 1}%
\special{ar 1350 1508 10 10 0  6.28318530717959E+0000}%
%
\special{pn 8}%
\special{ar 3060 2044 2220 860  4.1982242 4.6378167}%
%
\special{pn 8}%
\special{ar 3060 2044 2212 852  3.8259676 4.1447074}%
%
\special{pn 8}%
\special{ar 2320 1730 890 702  4.9145015 6.0690632}%
%
\special{pn 13}%
\special{pa 2896 1192}%
\special{pa 2642 1012}%
\special{fp}%
\special{sh 1}%
\special{pa 2642 1012}%
\special{pa 2684 1066}%
\special{pa 2684 1042}%
\special{pa 2708 1034}%
\special{pa 2642 1012}%
\special{fp}%
%
\special{pn 13}%
\special{pa 1342 1498}%
\special{pa 1088 1316}%
\special{fp}%
\special{sh 1}%
\special{pa 1088 1316}%
\special{pa 1130 1372}%
\special{pa 1130 1348}%
\special{pa 1154 1340}%
\special{pa 1088 1316}%
\special{fp}%
%
\special{pn 13}%
\special{pa 1498 1382}%
\special{pa 1244 1200}%
\special{fp}%
\special{sh 1}%
\special{pa 1244 1200}%
\special{pa 1288 1254}%
\special{pa 1288 1230}%
\special{pa 1310 1222}%
\special{pa 1244 1200}%
\special{fp}%
%
\special{pn 8}%
\special{pa 1356 1508}%
\special{pa 1490 1400}%
\special{fp}%
\special{sh 1}%
\special{pa 1490 1400}%
\special{pa 1426 1426}%
\special{pa 1448 1432}%
\special{pa 1452 1456}%
\special{pa 1490 1400}%
\special{fp}%
%
\special{pn 8}%
\special{pa 2896 1184}%
\special{pa 3030 1078}%
\special{fp}%
\special{sh 1}%
\special{pa 3030 1078}%
\special{pa 2966 1104}%
\special{pa 2988 1110}%
\special{pa 2990 1134}%
\special{pa 3030 1078}%
\special{fp}%
%
\special{pn 8}%
\special{ar 1588 2284 2376 422  4.0880892 5.5055163}%
%
\special{pn 8}%
\special{ar 1654 2210 1466 1578  3.3098433 4.3332690}%
%
\special{pn 8}%
\special{ar 1774 1598 3594 876  4.5243011 5.2873053}%
%
\special{pn 8}%
\special{ar 4040 2380 822 1662  3.3702437 4.2783353}%
\put(30.8100,-10.8500){\makebox(0,0)[lb]{$\xi_y$}}%
%
\special{pn 8}%
\special{pa 1654 2052}%
\special{pa 1432 1450}%
\special{dt 0.045}%
\special{sh 1}%
\special{pa 1432 1450}%
\special{pa 1436 1518}%
\special{pa 1450 1500}%
\special{pa 1474 1506}%
\special{pa 1432 1450}%
\special{fp}%
\put(16.2500,-20.9400){\makebox(0,0)[lt]{$\nu(y)$}}%
%
\special{pn 8}%
\special{pa 1490 598}%
\special{pa 1372 1292}%
\special{dt 0.045}%
\special{sh 1}%
\special{pa 1372 1292}%
\special{pa 1402 1230}%
\special{pa 1380 1240}%
\special{pa 1364 1224}%
\special{pa 1372 1292}%
\special{fp}%
%
\special{pn 8}%
\special{pa 1490 590}%
\special{pa 1162 1358}%
\special{dt 0.045}%
\special{sh 1}%
\special{pa 1162 1358}%
\special{pa 1208 1306}%
\special{pa 1184 1310}%
\special{pa 1170 1290}%
\special{pa 1162 1358}%
\special{fp}%
\put(13.6300,-5.4800){\makebox(0,0)[lb]{$({\rm grad}\,F)_{\nu(y)}$}}%
%
\special{pn 8}%
\special{pa 2940 606}%
\special{pa 2738 1078}%
\special{dt 0.045}%
\special{sh 1}%
\special{pa 2738 1078}%
\special{pa 2784 1024}%
\special{pa 2760 1028}%
\special{pa 2746 1008}%
\special{pa 2738 1078}%
\special{fp}%
\put(25.7300,-5.6400){\makebox(0,0)[lb]{$\tau_{f(y)}(({\rm grad}\,F)_{\nu(y)})\,\,(\in f_{\ast}(T_yM))$}}%
\put(37.2000,-12.3000){\makebox(0,0)[lb]{$f(M)$}}%
%
\special{pn 20}%
\special{sh 1}%
\special{ar 2888 1192 10 10 0  6.28318530717959E+0000}%
\special{sh 1}%
\special{ar 2888 1192 10 10 0  6.28318530717959E+0000}%
\put(29.1000,-12.5000){\makebox(0,0)[rt]{$f(y)$}}%
\put(21.1100,-13.1600){\makebox(0,0)[lt]{$\gamma_{f(y)}$}}%
\put(13.1800,-15.4800){\makebox(0,0)[lt]{$p_0$}}%
\put(4.6000,-11.7000){\makebox(0,0)[rb]{$T_{p_0}\widetilde M$}}%
%
\special{pn 8}%
\special{pa 924 590}%
\special{pa 1102 920}%
\special{dt 0.045}%
\special{sh 1}%
\special{pa 1102 920}%
\special{pa 1088 852}%
\special{pa 1078 872}%
\special{pa 1054 870}%
\special{pa 1102 920}%
\special{fp}%
\put(9.1600,-5.3900){\makebox(0,0)[rb]{$\widetilde M$}}%
%
\special{pn 8}%
\special{ar 1352 1510 212 164  0.0000000 6.2831853}%
%
\special{pn 8}%
\special{pa 3690 1190}%
\special{pa 3150 1450}%
\special{dt 0.045}%
\special{sh 1}%
\special{pa 3150 1450}%
\special{pa 3220 1440}%
\special{pa 3198 1428}%
\special{pa 3202 1404}%
\special{pa 3150 1450}%
\special{fp}%
%
\special{pn 8}%
\special{pa 570 880}%
\special{pa 1140 1490}%
\special{dt 0.045}%
\special{sh 1}%
\special{pa 1140 1490}%
\special{pa 1110 1428}%
\special{pa 1104 1452}%
\special{pa 1080 1456}%
\special{pa 1140 1490}%
\special{fp}%
\put(5.3000,-8.7000){\makebox(0,0)[rb]{$S^n(1)$}}%
%
\special{pn 8}%
\special{pa 480 1130}%
\special{pa 960 1640}%
\special{dt 0.045}%
\special{sh 1}%
\special{pa 960 1640}%
\special{pa 930 1578}%
\special{pa 924 1602}%
\special{pa 900 1606}%
\special{pa 960 1640}%
\special{fp}%
\end{picture}%
\hspace{5truecm}}

\vspace{0.3truecm}

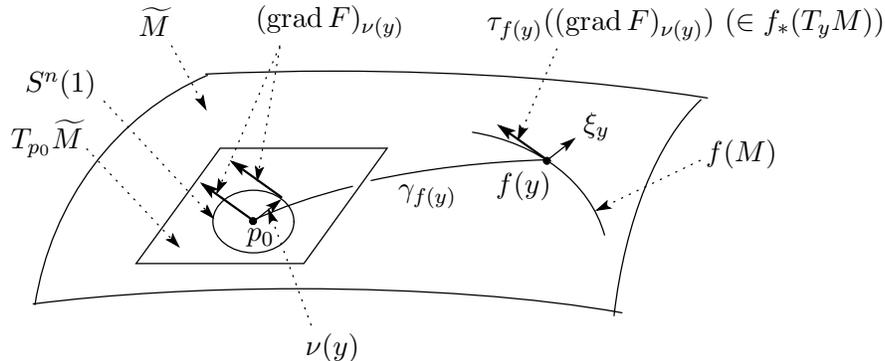
\centerline{{\bf Figure 1.}}

\newpage


\noindent
Take $v\in T_p\widetilde M$ 
and a curve $c:[0,\varepsilon)\to\widetilde M$ with $c'(0)=v$, we define 
a linear map $\tau^{\rm hol}_v:T_{p_0}\widetilde M\to T_{p_0}\widetilde M$ by 
$$\tau^{\rm hol}_v:=\left.\frac{d}{ds}\right\vert_{s=0}
\left(\tau_{\gamma_{c(0)}\cdot c\vert_{[0,s]}\cdot\gamma_{c(s)}^{-1}}\right),$$
where $\tau_{\gamma_{c(0)}\cdot c\vert_{[0,s]}\cdot\gamma_{c(s)}^{-1}}$ is 
the parallel translation along the product 
$\gamma_{c(0)}\cdot c\vert_{[0,s]}\cdot\gamma_{c(s)}^{-1}$
of the curves $\gamma_{c(0)},\,\,c\vert_{[0,s]}$ and $\gamma_{c(s)}^{-1}$.  
This linear map $\tau^{\rm hol}_v$ is described explicitly in the case where $\widetilde M$ is 
a symmetric space (see Lemma 2.1).  

Let $f_t\,\,(-\varepsilon<t<\varepsilon)$ be a ($C^{\infty}$-)variation of $f$ in 
${\rm Imm}(M,\widetilde M)$ and define 
$\widetilde f:M\times(-\varepsilon,\varepsilon)\to\widetilde M$ by 
$\widetilde f(x,t):=f_t(x)\,\,((x,t)\in M\times(-\varepsilon,\varepsilon))$.  
Denote by $g_t$ the metric on $M$ induced from $\langle\,\,,\,\,\rangle$ by $f_t$ 
and $\nabla^t$ the Riemannian connection of $g_t$.  Set $\nabla:=\nabla^0$.  
In the sequel, denote by $\langle\,\,,\,\,\rangle$ all metrics 
unless necessary.  Let $\pi_M$ be the natural projection of 
$M\times(-\varepsilon,\varepsilon)$ onto $M$.  
Denote by ${\widetilde{\nabla}}^f$ and 
${\widetilde{\nabla}}^{\widetilde f}$ the covariant derivative along 
$f$ and $\widetilde f$ for $\widetilde{\nabla}$, respectively.  
For a vector bundle $E$, denote by $\Gamma(E)$ the space of all ($C^{\infty}$-)sections 
of $E$.  
Also, denote by $TM$ the tangent bundle of $M$ and $\pi_M^{\ast}(TM)$ 
the bundle induced from $TM$ by $\pi_M$.  
For $X\in\Gamma(TM)$, we define $\overline X\in\Gamma(\pi_M^{\ast}(TM))$ by 
$\overline X_{(x,t)}:=X_x\,\,((x,t)\in M\times(-\varepsilon,\varepsilon))$.  
Let $\nabla'$ be the connection of $\pi_M^{\ast}(TM)$ satisfying 
$$(\nabla'_{\overline X}\overline Y)_{(x,t)}:=(\nabla^t_XY)_x\,\,\,
((x,t)\in M\times(-\varepsilon,\varepsilon))\quad{\rm and}\quad
\nabla'_{\frac{\partial}{\partial t}}\overline Y=0$$
for any $X,Y\in\Gamma(TM)$.  
Set 
$$\begin{array}{l}
\displaystyle{{\rm Imm}_{\rm b.f.}(f):=\{\hat f\in{\rm Imm}(M,\widetilde M)\,\vert\,
\exists\,{\rm a}\,\,{\rm boundary}\,\,{\rm fixing}\,\,{\rm variation}\,\,
f_t\,(0\leq t\leq1)}\\
\hspace{9.45truecm}\displaystyle{{\rm s.t.}\,\,f_0=f\,\,{\rm and}\,\,f_1
=\hat f\}}
\end{array}$$
and 
$$\begin{array}{l}
\displaystyle{{\rm Imm}_{\rm b.f.;v.p.}(f):=\{\hat f\in{\rm Imm}(M,\widetilde M)\,
\vert\,\exists\,{\rm a}\,\,{\rm boundary}\,\,{\rm fixing}\,\,{\rm and}\,\,
{\rm volume}\,\,{\rm preserving}}\\
\hspace{6.65truecm}\displaystyle{{\rm variation}\,\,f_t\,(0\leq t\leq1)\,\,
{\rm s.t.}\,\,f_0=f\,\,{\rm and}\,\,f_1=\hat f\}.}
\end{array}$$

We obtain the following first variational formula for ${\cal F}$.  

\vspace{0.3truecm}

\noindent
{\bf Theorem 1.1.} {\sl Let $f_t\,\,(-\varepsilon<t<\varepsilon)$ be a 
boundary-fixing variation of $f$ and $V=f_{\ast}(V_T)+\psi\xi$ the variational 
vector field of $f_t$, where $V_T\in\Gamma(TM)$ and $\psi\in C^{\infty}(M)$.  
Then we have 
$$\left.\frac{d}{dt}\right\vert_{t=0}{\cal F}(f_t)
=-\int_M\psi H_F\,dV,$$
where $H_F$ is the anisotropic mean curvature of $f$ and $dV$ is 
the volume element of the induced metric on $M$ by $f$.}

\vspace{0.3truecm}


\noindent
{\it Proof.} 
For simplicity, denote by $\nu_t$ the Gauss map of $f_t$.  
Denote by $g_t,\,A^t,\,h_t, H_t$ and $\xi_t$ the induced metric, the shape operator, the second 
fundamental form, the mean curvature and a unit normal vector field of $f_t$, 
respectively.  
Also, denote by $dV_t$ the volume element of $g_t$.  
For simplicity, set $A:=A^0,h:=h_0,H:=H_0$ and $\xi:=\xi_0$.  
Then we have 
$$\begin{array}{l}
\hspace{0.5truecm}
\displaystyle{\left.\frac{d}{dt}\right\vert_{t=0}{\cal F}(f_t)
=\int_{x\in M}\left.\frac{d}{dt}\right\vert_{t=0}(F(\nu_t(x))(dV_t)_x)}\\
\displaystyle{=\int_{x\in M}\left(\left(
\left.\frac{d}{dt}\right\vert_{t=0}F(\nu_t(x))\right)dV_x
+F(\nu(x))\left.\frac{d}{dt}\right\vert_{t=0}(dV_t)_x\right).}
\end{array}\leqno{(1.2)}$$
Also, we have 
$$\left.\frac{d}{dt}\right\vert_{t=0}(dV_t)_x
=(-H_x\psi(x)+({\rm div}\,V_T)_x)dV_x\leqno{(1.3)}$$
and 
$$\begin{array}{l}
\hspace{0.5truecm}
\displaystyle{\left.\frac{d}{dt}\right\vert_{t=0}F(\nu_t(x))
=\left\langle({\rm grad}\,F)_{\nu(x)},\,
\left.\frac{d}{dt}\right\vert_{t=0}\nu_t(x)\right\rangle.}
\end{array}\leqno{(1.4)}$$
Set $\beta(t):=f_t(x)$.  
By a simple calculation, we can show 
$$\widetilde{\nabla}^{\widetilde f}_{\frac{\partial}{\partial t}\vert_{t=0}}\xi_t
=-f_{\ast}\left(A(V_T)_x+({\rm grad}\,\psi)_x\right).
\leqno{(1.5)}$$
On the other hand, we have 
$$\begin{array}{l}
\displaystyle{\left.\frac{d}{dt}\right\vert_{t=0}\nu_t(x)
=\left.\frac{d}{dt}\right\vert_{t=0}\tau_{f_t(x)}^{-1}((\xi_t)_x)
}\\
\hspace{2.2truecm}\displaystyle{=\left.\frac{d}{dt}\right\vert_{t=0}
\tau_{\gamma_{f(x)}\cdot\beta\vert_{[0,t]}\cdot\gamma_{f_t(x)}^{-1}}
(\tau_{\gamma_{f(x)}\cdot\beta\vert_{[0,t]}}^{-1}((\xi_t)_x))}\\
\hspace{2.2truecm}\displaystyle{=\left(\left.\frac{d}{dt}\right\vert_{t=0}
\tau_{\gamma_{f(x)}\cdot\beta\vert_{[0,t]}\cdot\gamma_{f_t(x)}^{-1}}\right)
(\tau_{f(x)}^{-1}(\xi_x))}\\
\hspace{2.7truecm}\displaystyle{+\tau_{f(x)}^{-1}\left(\left.\frac{d}{dt}
\right\vert_{t=0}\tau_{\beta\vert_{[0,t]}^{-1}}((\xi_t)_x)\right)}\\
\hspace{2.2truecm}\displaystyle{=
\tau^{\rm hol}_{V_x}(\tau_{f(x)}^{-1}(\xi_x))
+\tau_{f(x)}^{-1}\left(
{\widetilde{\nabla}}^{\widetilde f}_{\frac{\partial}{\partial t}\vert_{t=0}}\xi_t\right).}
\end{array}\leqno{(1.6)}$$
Take a curve $c:[0,\varepsilon)\to M$ with $c'(0)=(V_T)_x$.  
Then we have 
$$\begin{array}{l}
\displaystyle{\tau_{f(x)}^{-1}(f_{\ast}A(V_T)_x)
=-\tau_{f(x)}^{-1}({\widetilde{\nabla}}^f_{(V_T)_x}\xi)}\\
\displaystyle{=-\left.\frac{d}{ds}\right\vert_{s=0}\tau_{\gamma_{(f\circ c)(s)}\cdot 
((f\circ c)\vert_{[0,s]})^{-1}\cdot\gamma_{f(x)}^{-1}}
(\tau_{(f\circ c)(s)}^{-1}(\xi_{c(s)}))}\\
\displaystyle{=\tau^{\rm hol}_{f_{\ast}((V_T)_x)}(\tau_{f(x)}^{-1}(\xi_x))
-\nu_{\ast x}((V_T)_x).}
\end{array}\leqno{(1.7)}$$
From $(1.4)-(1.7)$, we have 
$$\begin{array}{l}
\hspace{0.5truecm}\displaystyle{
\left.\frac{d}{dt}\right\vert_{t=0}F(\nu_t(x))}\\
\displaystyle{
=-({\rm div}(\psi(f_{\ast}^{-1}(\tau_{f(\cdot)}
({\rm grad}\,F)_{\nu(\cdot)})))_x+\psi(x)
({\rm div}(f_{\ast}^{-1}(\tau_{f(\cdot)}({\rm grad}\,F)_{\nu(\cdot)})))_x}\\
\hspace{0.5truecm}\displaystyle{
+\langle({\rm grad}\,F)_{\nu(x)},\tau^{\rm hol}_{V_x}
(\tau_{f(x)}^{-1}(\xi_x))\rangle
-\langle({\rm grad}\,F)_{\nu(x)},\tau^{\rm hol}_{f_{\ast}((V_T)_x)}
(\tau_{f(x)}^{-1}(\xi_x))\rangle}\\
\hspace{0.5truecm}\displaystyle{+({\rm div}((F\circ\nu)V_T))_x
-(F\circ\nu)(x)({\rm div}(V_T))_x,}
\end{array}\leqno{(1.8)}$$
where we use ${\rm div}(\phi Y)=\phi{\rm div}\,Y+\langle{\rm grad}\,\phi,\,Y
\rangle\,\,\,(\forall\,\phi\in C^{\infty}(M),\,\,\forall\,Y\in\Gamma(TM))$.  
Since $F$ is $\Phi$-invariant, we have 
$$\langle({\rm grad}\,F)_{\nu(x)},\tau^{\rm hol}_{V_x}(\tau_{f(x)}^{-1}(\xi_x))\rangle
=\langle({\rm grad}\,F)_{\nu(x)},\tau^{\rm hol}_{f_{\ast}((V_T)_x)}
(\tau_{f(x)}^{-1}(\xi_x))\rangle=0,\leqno{(1.9)}$$
where we note that $\tau^{\rm hol}_{V_x}$ and 
$\tau^{\rm hol}_{f_{\ast}((V_T)_x)}$ are elements of the holonomy algebra of $\widetilde M$ at $p_0$.  
Also, since $V$ vanishes on $\partial M$, we have 
$$\int_M{\rm div}(\psi(f_{\ast}^{-1}(\tau_{f(\cdot)}({\rm grad}\,F)_{\nu(\cdot)})))dV
=\int_M{\rm div}((F\circ\nu)V_T)dV=0.\leqno{(1.10)}$$
From $(1.2),\,(1.3),\,(1.8),\,(1.9)$ and $(1.10)$, 
we obtain the desried variational formula.\hspace{0.15truecm}q.e.d.

\vspace{0.5truecm}

\noindent
{\it Remark 1.1.} 
This first variational formula cannot be derived directly from the first variational formula by 
Lira and Melo (see the proof of [LM, Theorem 1]) for a (not necessarily holonomy invarinat) elliptic 
anisotropic surface energy.  

\vspace{0.5truecm}

From this first variational formula, we obtain the following result.  

\vspace{0.5truecm}

\noindent
{\bf Theorem 1.2.} 
{\sl {\rm(i)} $f$ is a critical point of 
${\cal F}\vert_{{\rm Imm}_{\rm b.f.}(f)}$ if and only if $H_F$ vanishes.  

{\rm(ii)} $f$ is a critical point of 
${\cal F}\vert_{{\rm Imm}_{\rm b.f.;v.p.}(f)}$ if and only if $H_F$ is 
constant.}

\vspace{0.5truecm}

\noindent
{\it Proof.} The statement (i) follows from the first variational formula in Theorem 1.1 directly.  
When the variation $f_t$ ($-\varepsilon<t<\varepsilon$) is volume-preserving, 
we have $\int_M\psi dV=0$, where $\psi$ is as in Theorem 1.1.  
Hence the statement (ii) also follows from the first variational formula in Theorem 1.1 directly.  
\hspace{10.2truecm}q.e.d.

\section{Anisotropic equifocal hypersurfaces 
and anisotropic isoparametric hypersurfaces} 
We use the notations in the previous section.  
In this section, we shall introduce the notions of an anisotropic convex hypersurface, 
an anisotropic equifocal hypersurface and an anisotropic isoparametric hypersurface 
for the holonomy invariant elliptic anisotropic surface energy ${\cal F}$.  
Assume that $\partial\,M=\emptyset$.  
Define a transversal vector field $\xi_F$ of $f$ by 
$$(\xi_F)_x:=(F\circ\nu)(x)\xi_x+\tau_{f(x)}({\rm grad}\,F)_{\nu(x)}
\quad\,\,\,\,(x\in M).\leqno{(2.1)}$$
We call $\xi_F$ a {\it anisotropic transversal vector field} of $f$.  
Take $X\in T_xM$.  Let $\beta:[0,\varepsilon)\to M$ be a curve with 
$\beta'(0)=X$.  Then, since $F$ is $\Phi$-invariant, we have 
$$\begin{array}{l}
\displaystyle{{\widetilde{\nabla}}^f_X\xi_F=X(F\circ\nu)\xi_x
-(F\circ\nu)(x)f_{\ast}(A_xX)+\tau_{f(x)}((\tau^{\rm hol}_X)^{-1}
(({\rm grad}\,F)_{\nu(x)}))}\\
\hspace{1.7truecm}\displaystyle{+\tau_{f(x)}\left(
\nabla^S_{\nu_{\ast}(X)}{\rm grad}\,F-\langle\nu_{\ast}(X),{\rm grad}\,F
\rangle\nu(x)\right)}\\
\hspace{1.16truecm}\displaystyle{=-(F\circ\nu)(x)f_{\ast}(A_xX)
+f_{\ast}\left(\nabla_X(f_{\ast}^{-1}((\tau_{f(\cdot)}({\rm grad}\,F)_{\nu(\cdot)})))\right).}
\end{array}\leqno{(2.2)}$$
So, define a $(1,1)$-tensor field $A^F$ on $M$ by 
$$A^F_xX=(F\circ\nu)(x)A_xX-\nabla_X(f_{\ast}^{-1}(\tau_{f(\cdot)}
({\rm grad}\,F)_{\nu(\cdot)}))\leqno{(2.3)}$$
for any $x\in M$ and $X\in T_xM$.  
We call $A^F$ a {\it anisotropic shape operator} of $f$ and 
the eigenvalues of $A^F_x$ {\it anisotropic principal curvatures of} $f$ 
at $x$.  It is easy to show that ${\rm Tr}\,A^F$ coincides with $H_F$.  
If all {\it anisotropic principal curvatures of} $f$ are positive 
(or negative) at each point of $M$, then we say that 
$f:M\hookrightarrow\widetilde M$ is {\it anisotropic convex}.  
Note that, in the case where $\widetilde M$ is a $(n+1)$-dimensional Euclidean space 
${\Bbb R}^{n+1}$, we have 
$$f_{\ast}(\nabla_X(f_{\ast}^{-1}(\tau_{f(\cdot)}({\rm grad}\,F)_{\nu(\cdot)})))
=\nabla^S_{\nu_{\ast}(X)}{\rm grad}\,F=-\nabla^S_{f_{\ast}(A_xX)}{\rm grad}\,F$$
and hence 
$$f_{\ast}(A^F_xX)=(F\circ\nu)(x)f_{\ast}(A_xX)
+\nabla^S_{f_{\ast}(A_xX)}{\rm grad}\,F,$$
where we identify $f_{\ast}(T_xM)$ with $T_{\nu(x)}S^n(1)$ under the identification of 
$T_{f(x)}{\Bbb R}^{n+1}$ and $T_{p_0}{\Bbb R}^{n+1}$.  
Denote by $R$ the curvature tensor of $\widetilde M$ and 
$R(\xi_F)$ the anisotropic normal Jacobi operator $R(\cdot,\xi_F)\xi_F$.  
Let $\gamma^F_x$ be the geodesic in $\widetilde M$ 
whose initial velocity vector is equal to $(\xi_F)_x$.  
Take $X\in T_xM$.  
Let $Y_X$ be the Jacobi field along $\gamma^F_x$ with $Y_X(0)=f_{\ast}(X)$ and 
$Y_X'(0)=-f_{\ast}(A^F_xX)$.  We call $Y_X$ an {\it anisotropic} $M$-{\it Jacobi field}.  
If $Y_X(s_0)=0$ for some $X(\not=0)\in T_xM$, then we call $s_0$ (resp. $\gamma^F_x(s_0)$) 
an {\it anisotropic focal radius} (resp. an {\it anisotropic focal point}) of $f$ at $x$.  
Also, for a anisotropic focal radius $s_0$ of $f$ at $x$, we call 
$\{X\in T_xM\,\vert\,Y_X(s_0)=0\}$ the {\it nullity space} of $s_0$ and 
its dimension the {\it multiplicity} of $s_0$.  
If the set of all anisotropic focal radii of $f$ at $x$ is 
independent of the choice of $x\in M$, then we call $f:M\hookrightarrow\widetilde M$ 
an {\it anisotropic equifocal hypersurface}.  
Also, if the set of all anisotropic principal curvatures of $f$ at 
$x$ is independent of the choice of $x\in M$, then we call 
$f:M\hookrightarrow\widetilde M$ a {\it hypersurface with constant anisotropic principal curvatures}.  
Next we shall introduce the notion of an anisotropic isoparametric hypersurface.  
Define a map $f_t:M\to\widetilde M$ by 
$$f_t(x):=\exp_{f(x)}(t(\xi_F)_x)\quad\,\,(x\in M),$$
where $\exp_{f(x)}$ is the exponential map of $\widetilde M$ at $f(x)$.  
It is easy to show that $f_t$ is an immersion for each $t$ sufficiently close 
to zero.  If $f_t$ is an immersion, then we call 
$f_t:M\hookrightarrow\widetilde M$ an {\it anisotropic parallel hypersurface} of 
$f:M\hookrightarrow\widetilde M$ {\it of distance} $t$.  
Also, if $f_t$ is not an immersion but the differential of $f_t$ at each point of $M$ is of 
constant rank, then we call $f_t(M)$ an {\it anisotropic focal submanifold of} 
$f:M\hookrightarrow\widetilde M$ ({\it corresponding to anisotropic focal radius} $t$).  
If, for each $t$ sufficiently close to zero, $f_t:M\hookrightarrow\widetilde M$ is 
of constant anisotropic mean curvature, we call 
$f:M\hookrightarrow\widetilde M$ an {\it anisotropic isoparametric hypersurface}.  

We consider the case where $\widetilde M$ is a symmetric space.  
Then a Jacobi field $Y$ along a geodesic $\gamma$ in $\widetilde M$ is described as 
$$Y(s)=\tau_{\gamma\vert_{[0,s]}}\left(D^{co}_{s\gamma'(0)}(Y(0))
+sD^{si}_{s\gamma'(0)}(Y'(0))\right),\leqno{(2.4)}$$
where $D^{co}_{s\gamma'(0)}$ (resp. $D^{si}_{s\gamma'(0)}$) is given by 
$$
D^{co}_{s\gamma'(0)}:=\cos(s\sqrt{R(\gamma'(0))})\quad\,
\left({\rm resp.}\,\,
D^{si}_{s\gamma'(0)}:=\frac{\sin(s\sqrt{R(\gamma'(0))})}
{s\sqrt{R(\gamma'(0))}}\right).
$$
In particular, the anisotropic $M$-Jacobi field $Y_X$ of $f$ is described as 
$$Y_X(s)=\tau_{\gamma^F_x\vert_{[0,s]}}\left(D^{co}_{s(\xi_F)_x}(f_{\ast}X)
-sD^{si}_{s(\xi_F)_x}(f_{\ast}(A^F_xX))\right),\leqno{(2.5)}$$
where $D^{co}_{s(\xi_F)_x}$ (resp. $D^{si}_{s(\xi_F)_x}$) is defined in a similar way to 
$D^{co}_{s\gamma'(0)}$ (resp. $D^{si}_{s\gamma'(0)}$).  
According to $(2.5)$, the anisotropic focal radii of $f$ at $x$ coincide with zero points of 
the function 
$$\rho(s):={\rm det}\left(D^{co}_{s(\xi_F)_x}\circ f_{\ast}
-s(D^{si}_{s(\xi_F)_x}\circ f_{\ast}\circ A^F_x)\right).$$
In particular, in the case where $\widetilde M$ is a Euclidean space, 
$D^{co}_{s(\xi_F)_x}=D^{si}_{s(\xi_F)_x}={\rm id}$ and hence 
the anisotropic focal radius of $f$ at $x$ are equal to the inverse numbers of 
anisotropic principal curvatures of $f$ at $x$.  

At the end of this section, we give an explicit description of the linear map $\tau^{\rm hol}_v$ 
defined in the previous secition in the case where $\widetilde M$ is a symmetric space.  

\vspace{0.5truecm}

\noindent
{\bf Lemma 2.1.} {\sl For $w\in T_{p_0}\widetilde M$, we have 
$$\tau^{hol}_v(w)=R_{p_0}\left(\gamma'_p(0),
\frac{{\rm id}-D^{co}_{\gamma'_p(0)}}{{\rm ad}(\gamma'_p(0))^2}
\left((\exp_{p_0})_{\ast\gamma'_p(0)}^{-1}(v)\right)\right)w.$$
}

\vspace{0.5truecm}

\noindent
{\it Proof.} 
Set $\overline v:=(\exp_{p_0})_{\ast\gamma'_p(0)}^{-1}(v)$.  
Define a $2$-parameter map $\delta:[0,1]^2\to\widetilde M$ by 
$$\delta(s,t):=\exp_{p_0}\left(s(\gamma'_p(0)+t\overline v)\right)\quad((s,t)\in[0,1]^2).$$
Let $Y$ be the vector field along $\gamma_p$ defined by 
$\displaystyle{Y:=\left.\frac{\partial\delta}{\partial t}\right\vert_{t=0}}$.  
Since $Y$ is the Jacobi field along $\gamma_p$ with $Y(0)=0$ and $Y'(0)=\overline v$, 
it is described as 
$$Y(s)=\tau_{\gamma_p\vert_{[0,s]}}\left(sD^{si}_{s\gamma'_p(0)}\overline v)\right).
\leqno{(2.6)}$$
Let $\widetilde w$ be the vector field along $\delta$ with $\widetilde w_{(0,0)}=w$ such that 
$s\mapsto\widetilde w_{(s,t)}$ is parallel along $s\mapsto\delta(s,t)$ for each $t\in[0,1]$.  
Define a $3$-parameter map $\widetilde{\delta}:[0,1]^3\to\widetilde M$ by 
$$\widetilde{\delta}(s,t,u):=\exp_{\delta(s,t)}\left(u\widetilde w_{(s,t)}\right)
\quad((s,t,u)\in[0,1]^3).$$
Clearly we have 
$$\tau^{hol}_v(w)=-\tau_p^{-1}\left(
\widetilde{\nabla}^{\widetilde{\delta}}_{\frac{\partial}{\partial t}}
\left.\frac{\partial\widetilde{\delta}}{\partial u}
\right\vert_{s=1,t=0,u=0}\right).\leqno{(2.7)}$$
Define a vector field $Z$ along $\delta$ by 
$$Z:=\widetilde{\nabla}^{\widetilde{\delta}}_{\frac{\partial}{\partial t}}
\left.\frac{\partial\widetilde{\delta}}{\partial u}\right\vert_{u=0}.$$
Then, by using $(2.6)$, we can show 
$$\begin{array}{l}
\displaystyle{\left.\left(\widetilde{\nabla}^{\delta}_{\frac{\partial}{\partial s}}Z\right)
\right\vert_{t=0}=\left.\left(\widetilde{\nabla}^{\delta}_{\frac{\partial}{\partial t}}
\widetilde{\nabla}^{\delta}_{\frac{\partial}{\partial s}}\widetilde w\right)\right\vert_{t=0}
+R_{p_0}\left(\left.\frac{\partial\delta}{\partial s}\right\vert_{t=u=0},
\left.\frac{\partial\delta}{\partial t}\right\vert_{t=u=0}\right)\widetilde w\vert_{t=0}}\\
\displaystyle{=R_{p_0}(\gamma_p'(s),Y(s))\widetilde w\vert_{t=0}
=\tau_{\gamma_p\vert_{[0,s]}}\left(
R_{p_0}\left(\gamma_p'(0),sD^{si}_{s\gamma_p'(0)}\overline v\right)w\right).}
\end{array}$$
Also, we have $Z\vert_{s=t=0}=0$.  
Hence we obtain 
$$Z_{(s,t)}=\tau_{\gamma_p\vert_{[0,s]}}\left(\int_0^sR_{p_0}(\gamma_p'(0),sD^{si}_{s\gamma_p'(0)}
\overline v)w\,ds\right).$$
In particular, we obtain 
$$\begin{array}{l}
\displaystyle{Z_{(1,0)}=\tau_p\left(\int_0^1R_{p_0}(\gamma_p'(0),
sD^{si}_{s\gamma_p'(0)}\overline v)w\,ds\right)}\\
\displaystyle{=\tau_p\left(R_{p_0}\left(\gamma_p'(0),
\left(\int_0^1sD^{si}_{s\gamma_p'(0)}ds\right)(\overline v)\right)w\right)}\\
\displaystyle{=\tau_p\left(R_{p_0}\left(\gamma_p'(0),
\frac{D^{co}_{\gamma_p'(0)}-{\rm id}}{{\rm ad}(\gamma_p'(0))^2}(\overline v)\right)w\right)}
\end{array}.$$
This relation together with $(2.7)$ implies the desired relation.  \hspace{3.5truecm}q.e.d.

\section{Anisotropic tubes} 
We use the notations in Sections 1 and 2.  
In this section, we introduce the notion of anisotropic tube over a submanifold.  
Let $B$ be an embedded submanifold in $\widetilde M$ and $\pi_B:T^{\perp_1}B\to B$ 
the unit normal bundle of $B$.  
For a positive number $r$, we define $f^F_{T,r}:T^{\perp_1}B\to\widetilde M$ by 
$$f^F_{B,r}(v):=\exp_{\pi_B(v)}\left(r\left(\widetilde F(v)v
+({\rm grad}(F_{\pi_B(v)}))_v\right)\right)\quad(v\in T^{\perp_1}B).$$
Set $t^F_r(B):=f^F_{B,r}(T^{\perp_1}B)$.  If $f^F_{B,r}$ is an immersion, then we call 
$t^F_r(B)$ the {\it anisotropic tube over} $B$ {\it of radius $r$}.  
For an anisotropic tube, we can show the following fact.  

\vspace{0.5truecm}

\noindent
{\bf Theorem 3.1.} {\sl Let $B$ be a complete embedded submanifold in $\widetilde M$ and 
$f:M\hookrightarrow\widetilde M$ a complete hypersurface in $\widetilde M$.  
If $B$ is an anisotropic focal submanifold of $f(M)$ corresponding to anisotropic focal radius $r$, 
then $f(M)=t_{-r}^F(B)$ holds.}

\vspace{0.5truecm}

\noindent
{\it Proof.} 
By the assumption, we have $B=f_r(M)$.  Set $\varepsilon:=\frac{r}{\vert r\vert}$.  
Let $\xi(s)$ and $\xi_F(s)$ ($0\leq\varepsilon s<\varepsilon r$) be the unit normal vector field 
and the anisotropic normal vector field of the parallel hypersurface $f_s(M)(=t_{s-r}^F(B))$ of $f(M)$, 
respectively.  
Take $x\in M$ and set $p:=f_r(x)$.  Since $f_s(x)=\gamma_{(\xi_F(0))_x}(s)$, we have 
$(\xi_F(s))_x=\gamma'_{(\xi_F(0))_x}(s)$.  
Also, since $(\xi_F(s))_x=\widetilde F(\xi(s)_x)\xi(s)_x
+({\rm grad}(F_{f_s(x)}))_{\xi(s)_x}$, 
we obtain 
$$\gamma'_{(\xi_F(0))_x}(r)=
\widetilde F(\lim_{s\to r}\xi(s)_x)\lim_{s\to r}\xi(s)_x
+({\rm grad}(F\vert_p))_{\lim\limits_{s\to r}\xi(s)_x}.$$
Set $v:=\lim\limits_{s\to r}\xi(s)_x$.  Then we have 
$$\begin{array}{l}
\displaystyle{f(x)=\gamma_{(\xi_F(0))_x}(0)=\gamma_{\gamma'_{(\xi_F(0))_x}(r)}(-r)
=\exp_p(-r\gamma'_{(\xi_F(0))_x}(r))}\\
\displaystyle{=\exp_p\left(-r\left(\widetilde F(v)v+({\rm grad}(F_p))_v\right)\right)
=f^F_{B,-r}(v)\in t^F_{-r}(B).}
\end{array}$$
Thus it follows from the arbitrariness of $x$ that $f(M)\subset t_{-r}^F(B)$.  
Furthermore, it follows from the completeness of $f(M)$ that $f(M)=t_{-r}^F(B)$.  
\hspace{4.95truecm}q.e.d.

\vspace{0.5truecm}

\centerline{
\unitlength 0.1in
\begin{picture}( 19.7500, 18.9500)( 19.0500,-21.6200)
%
\special{pn 8}%
\special{ar 2520 1920 384 120  6.2831853 6.2831853}%
\special{ar 2520 1920 384 120  0.0000000 3.1415927}%
%
\special{pn 8}%
\special{ar 2550 1920 346 120  4.8395485 6.2831853}%
%
\special{pn 8}%
\special{pa 2136 1910}%
\special{pa 2148 1882}%
\special{pa 2170 1860}%
\special{pa 2198 1844}%
\special{pa 2228 1830}%
\special{pa 2258 1820}%
\special{pa 2288 1812}%
\special{pa 2320 1806}%
\special{pa 2352 1802}%
\special{pa 2384 1798}%
\special{pa 2416 1794}%
\special{pa 2448 1794}%
\special{pa 2476 1792}%
\special{sp}%
%
\special{pn 8}%
\special{ar 3282 1780 376 1236  2.8263112 3.8891710}%
%
\special{pn 8}%
\special{ar 2514 1696 376 1236  2.8263112 3.8891710}%
%
\special{pn 8}%
\special{ar 2918 1734 374 1236  2.8276677 3.8880495}%
%
\special{pn 8}%
\special{pa 2944 1238}%
\special{pa 2928 1266}%
\special{pa 2904 1286}%
\special{pa 2876 1300}%
\special{pa 2846 1310}%
\special{pa 2814 1316}%
\special{pa 2782 1322}%
\special{pa 2750 1326}%
\special{pa 2718 1330}%
\special{pa 2686 1332}%
\special{pa 2654 1332}%
\special{pa 2622 1332}%
\special{pa 2590 1332}%
\special{pa 2560 1330}%
\special{pa 2528 1326}%
\special{pa 2496 1324}%
\special{pa 2464 1318}%
\special{pa 2432 1314}%
\special{pa 2400 1308}%
\special{pa 2370 1300}%
\special{pa 2338 1292}%
\special{pa 2308 1282}%
\special{pa 2278 1272}%
\special{pa 2250 1258}%
\special{pa 2222 1242}%
\special{pa 2196 1222}%
\special{pa 2180 1196}%
\special{pa 2178 1180}%
\special{sp}%
%
\special{pn 8}%
\special{pa 2620 1088}%
\special{pa 2652 1092}%
\special{pa 2682 1098}%
\special{pa 2714 1102}%
\special{pa 2746 1108}%
\special{pa 2776 1116}%
\special{pa 2808 1126}%
\special{pa 2838 1136}%
\special{pa 2866 1148}%
\special{pa 2894 1166}%
\special{pa 2918 1184}%
\special{pa 2940 1208}%
\special{pa 2946 1236}%
\special{sp}%
%
\special{pn 8}%
\special{pa 2178 1190}%
\special{pa 2190 1160}%
\special{pa 2214 1140}%
\special{pa 2242 1126}%
\special{pa 2270 1114}%
\special{pa 2302 1106}%
\special{pa 2332 1098}%
\special{pa 2364 1092}%
\special{pa 2396 1090}%
\special{pa 2428 1088}%
\special{pa 2460 1086}%
\special{pa 2492 1086}%
\special{pa 2522 1086}%
\special{sp}%
%
\special{pn 20}%
\special{sh 1}%
\special{ar 2560 1910 10 10 0  6.28318530717959E+0000}%
\special{sh 1}%
\special{ar 2560 1910 10 10 0  6.28318530717959E+0000}%
%
\special{pn 8}%
\special{ar 2544 1580 626 146  6.2831853 6.2831853}%
\special{ar 2544 1580 626 146  0.0000000 3.1415927}%
%
\special{pn 8}%
\special{ar 2544 1580 626 158  4.8397743 6.2831853}%
%
\special{pn 8}%
\special{ar 2544 1580 624 158  3.1415927 4.5853358}%
%
\special{pn 8}%
\special{pa 3176 2158}%
\special{pa 3174 2126}%
\special{pa 3170 2094}%
\special{pa 3166 2062}%
\special{pa 3164 2030}%
\special{pa 3162 1998}%
\special{pa 3158 1966}%
\special{pa 3156 1934}%
\special{pa 3156 1902}%
\special{pa 3154 1870}%
\special{pa 3154 1838}%
\special{pa 3152 1806}%
\special{pa 3152 1774}%
\special{pa 3152 1742}%
\special{pa 3152 1710}%
\special{pa 3152 1678}%
\special{pa 3154 1646}%
\special{pa 3154 1614}%
\special{pa 3156 1582}%
\special{pa 3156 1550}%
\special{pa 3158 1518}%
\special{pa 3160 1486}%
\special{pa 3164 1454}%
\special{pa 3164 1422}%
\special{pa 3168 1392}%
\special{pa 3172 1360}%
\special{pa 3174 1328}%
\special{pa 3178 1296}%
\special{pa 3184 1264}%
\special{pa 3188 1232}%
\special{pa 3192 1200}%
\special{pa 3198 1170}%
\special{pa 3202 1138}%
\special{pa 3208 1106}%
\special{pa 3214 1074}%
\special{pa 3220 1044}%
\special{pa 3226 1012}%
\special{pa 3232 980}%
\special{pa 3240 950}%
\special{pa 3248 918}%
\special{pa 3256 888}%
\special{pa 3266 856}%
\special{pa 3274 826}%
\special{pa 3282 796}%
\special{pa 3294 766}%
\special{pa 3304 734}%
\special{pa 3316 704}%
\special{pa 3326 674}%
\special{pa 3326 672}%
\special{sp}%
%
\special{pn 8}%
\special{pa 3250 910}%
\special{pa 3236 938}%
\special{pa 3210 958}%
\special{pa 3182 972}%
\special{pa 3152 984}%
\special{pa 3120 990}%
\special{pa 3090 998}%
\special{pa 3058 1004}%
\special{pa 3026 1010}%
\special{pa 2994 1014}%
\special{pa 2962 1016}%
\special{pa 2930 1020}%
\special{pa 2898 1022}%
\special{pa 2866 1024}%
\special{pa 2834 1024}%
\special{pa 2802 1024}%
\special{pa 2770 1024}%
\special{pa 2738 1024}%
\special{pa 2706 1024}%
\special{pa 2674 1024}%
\special{pa 2642 1024}%
\special{pa 2610 1022}%
\special{pa 2580 1020}%
\special{pa 2548 1018}%
\special{pa 2516 1014}%
\special{pa 2484 1012}%
\special{pa 2452 1008}%
\special{pa 2420 1004}%
\special{pa 2388 1000}%
\special{pa 2356 994}%
\special{pa 2326 988}%
\special{pa 2294 984}%
\special{pa 2262 976}%
\special{pa 2230 970}%
\special{pa 2200 962}%
\special{pa 2170 952}%
\special{pa 2138 942}%
\special{pa 2110 930}%
\special{pa 2080 918}%
\special{pa 2052 904}%
\special{pa 2026 884}%
\special{pa 2008 858}%
\special{pa 2004 840}%
\special{sp}%
%
\special{pn 8}%
\special{pa 2730 750}%
\special{pa 2762 754}%
\special{pa 2794 756}%
\special{pa 2826 760}%
\special{pa 2858 764}%
\special{pa 2890 768}%
\special{pa 2922 772}%
\special{pa 2952 778}%
\special{pa 2984 782}%
\special{pa 3014 792}%
\special{pa 3046 798}%
\special{pa 3078 804}%
\special{pa 3106 816}%
\special{pa 3138 826}%
\special{pa 3166 838}%
\special{pa 3194 852}%
\special{pa 3220 872}%
\special{pa 3238 896}%
\special{pa 3240 910}%
\special{sp}%
%
\special{pn 8}%
\special{pa 1990 850}%
\special{pa 2006 824}%
\special{pa 2032 806}%
\special{pa 2060 792}%
\special{pa 2090 784}%
\special{pa 2120 774}%
\special{pa 2152 768}%
\special{pa 2184 762}%
\special{pa 2214 756}%
\special{pa 2246 752}%
\special{pa 2278 748}%
\special{pa 2310 746}%
\special{pa 2342 744}%
\special{pa 2374 744}%
\special{pa 2406 742}%
\special{pa 2438 742}%
\special{pa 2470 740}%
\special{pa 2502 742}%
\special{pa 2534 740}%
\special{pa 2566 740}%
\special{pa 2598 742}%
\special{pa 2630 744}%
\special{pa 2630 744}%
\special{sp}%
%
\special{pn 8}%
\special{pa 1918 2092}%
\special{pa 1916 2060}%
\special{pa 1914 2028}%
\special{pa 1912 1996}%
\special{pa 1912 1964}%
\special{pa 1910 1932}%
\special{pa 1908 1900}%
\special{pa 1908 1868}%
\special{pa 1908 1836}%
\special{pa 1906 1804}%
\special{pa 1906 1772}%
\special{pa 1906 1740}%
\special{pa 1906 1708}%
\special{pa 1906 1676}%
\special{pa 1906 1644}%
\special{pa 1908 1612}%
\special{pa 1908 1580}%
\special{pa 1910 1548}%
\special{pa 1910 1516}%
\special{pa 1912 1486}%
\special{pa 1914 1454}%
\special{pa 1916 1422}%
\special{pa 1918 1390}%
\special{pa 1920 1358}%
\special{pa 1922 1326}%
\special{pa 1924 1294}%
\special{pa 1928 1262}%
\special{pa 1932 1230}%
\special{pa 1936 1198}%
\special{pa 1938 1166}%
\special{pa 1942 1134}%
\special{pa 1946 1102}%
\special{pa 1950 1070}%
\special{pa 1956 1040}%
\special{pa 1962 1008}%
\special{pa 1966 976}%
\special{pa 1970 944}%
\special{pa 1978 914}%
\special{pa 1982 882}%
\special{pa 1988 850}%
\special{pa 1994 818}%
\special{pa 2002 788}%
\special{pa 2008 756}%
\special{pa 2016 726}%
\special{pa 2024 694}%
\special{pa 2032 664}%
\special{pa 2040 632}%
\special{pa 2050 602}%
\special{pa 2058 572}%
\special{pa 2060 566}%
\special{sp}%
%
\special{pn 8}%
\special{pa 3484 842}%
\special{pa 2870 1368}%
\special{dt 0.045}%
\special{sh 1}%
\special{pa 2870 1368}%
\special{pa 2934 1340}%
\special{pa 2910 1334}%
\special{pa 2908 1310}%
\special{pa 2870 1368}%
\special{fp}%
%
\special{pn 8}%
\special{pa 3560 1120}%
\special{pa 3090 1332}%
\special{dt 0.045}%
\special{sh 1}%
\special{pa 3090 1332}%
\special{pa 3158 1322}%
\special{pa 3138 1310}%
\special{pa 3142 1286}%
\special{pa 3090 1332}%
\special{fp}%
\put(36.0800,-11.5600){\makebox(0,0)[lb]{$t^F_r(B)$ (outside tube)}}%
\put(35.2000,-8.5000){\makebox(0,0)[lb]{$\exp^{\perp}(T^{\perp_1}B)$ (inside tube)}}%
%
\special{pn 8}%
\special{ar 3282 1756 740 1688  3.6622686 3.8729035}%
%
\special{pn 8}%
\special{ar 2878 1718 740 1688  3.6617577 3.8729035}%
%
\special{pn 8}%
\special{ar 3646 1810 740 1688  3.6617577 3.8729035}%
%
\special{pn 20}%
\special{sh 1}%
\special{ar 2910 1920 10 10 0  6.28318530717959E+0000}%
\special{sh 1}%
\special{ar 2910 1920 10 10 0  6.28318530717959E+0000}%
%
\special{pn 20}%
\special{sh 1}%
\special{ar 3166 1590 10 10 0  6.28318530717959E+0000}%
\special{sh 1}%
\special{ar 3166 1590 10 10 0  6.28318530717959E+0000}%
%
\special{pn 8}%
\special{pa 3512 1368}%
\special{pa 3176 1590}%
\special{dt 0.045}%
\special{sh 1}%
\special{pa 3176 1590}%
\special{pa 3244 1570}%
\special{pa 3222 1560}%
\special{pa 3222 1536}%
\special{pa 3176 1590}%
\special{fp}%
%
\special{pn 8}%
\special{pa 3820 1860}%
\special{pa 2918 1916}%
\special{dt 0.045}%
\special{sh 1}%
\special{pa 2918 1916}%
\special{pa 2986 1932}%
\special{pa 2972 1912}%
\special{pa 2984 1892}%
\special{pa 2918 1916}%
\special{fp}%
\put(38.8000,-17.8000){\makebox(0,0)[lt]{$\exp^{\perp}\,v$}}%
\put(35.7900,-13.8600){\makebox(0,0)[lb]{$f^F_{B,r}(v)$}}%
%
\special{pn 8}%
\special{pa 2494 484}%
\special{pa 2696 714}%
\special{dt 0.045}%
\special{sh 1}%
\special{pa 2696 714}%
\special{pa 2668 650}%
\special{pa 2662 674}%
\special{pa 2638 676}%
\special{pa 2696 714}%
\special{fp}%
\put(25.1300,-4.3700){\makebox(0,0)[rb]{$B$}}%
%
\special{pn 8}%
\special{pa 2180 1180}%
\special{pa 2940 1230}%
\special{dt 0.045}%
%
\special{pn 8}%
\special{pa 2140 1900}%
\special{pa 2900 1910}%
\special{dt 0.045}%
%
\special{pn 8}%
\special{pa 2560 1910}%
\special{pa 3160 1600}%
\special{dt 0.045}%
%
\special{pn 8}%
\special{pa 2550 1900}%
\special{pa 1910 1610}%
\special{dt 0.045}%
%
\special{pn 8}%
\special{pa 2590 1200}%
\special{pa 3240 950}%
\special{dt 0.045}%
%
\special{pn 8}%
\special{pa 2580 1210}%
\special{pa 1990 860}%
\special{dt 0.045}%
%
\special{pn 20}%
\special{sh 1}%
\special{ar 2590 1210 10 10 0  6.28318530717959E+0000}%
\special{sh 1}%
\special{ar 2590 1210 10 10 0  6.28318530717959E+0000}%
%
\special{pn 20}%
\special{sh 1}%
\special{ar 3250 930 10 10 0  6.28318530717959E+0000}%
\special{sh 1}%
\special{ar 3250 930 10 10 0  6.28318530717959E+0000}%
%
\special{pn 20}%
\special{sh 1}%
\special{ar 1910 1600 10 10 0  6.28318530717959E+0000}%
\special{sh 1}%
\special{ar 1910 1600 10 10 0  6.28318530717959E+0000}%
%
\special{pn 20}%
\special{sh 1}%
\special{ar 1990 840 10 10 0  6.28318530717959E+0000}%
\special{sh 1}%
\special{ar 1990 840 10 10 0  6.28318530717959E+0000}%
%
\special{pn 20}%
\special{sh 1}%
\special{ar 2130 1900 10 10 0  6.28318530717959E+0000}%
\special{sh 1}%
\special{ar 2130 1900 10 10 0  6.28318530717959E+0000}%
%
\special{pn 20}%
\special{sh 1}%
\special{ar 2170 1180 10 10 0  6.28318530717959E+0000}%
\special{sh 1}%
\special{ar 2170 1180 10 10 0  6.28318530717959E+0000}%
%
\special{pn 20}%
\special{sh 1}%
\special{ar 2950 1230 10 10 0  6.28318530717959E+0000}%
\special{sh 1}%
\special{ar 2950 1230 10 10 0  6.28318530717959E+0000}%
\end{picture}%
\hspace{2truecm}}

\vspace{0.5truecm}

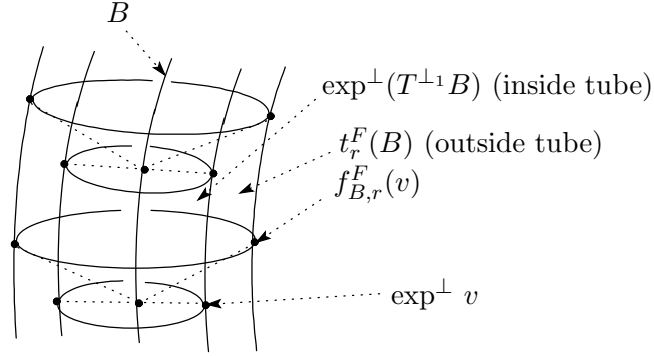
\centerline{{\bf Figure 2.}}

\section{The anisotropic geodesic spheres}
We use the notations in Sections 1-3.  
For simplicity, set ${\rm Exp}:=\exp_{p_0}$.  
In this section, we consider the anisotropic tube $t^F_r(p_0)$ over $\{p_0\}$.  
We call this tube the {\it anisotropic geodesic sphere of} $r$ {\it centered at} $p_0$.  
For example, in the setting of [KP1-3, Palm], the Wulff shape is the anisotropic geodesic sphere 
of radius $1$ centered at the origin.  
Define ${\widehat f}^F_{p_0,r}:S^n(1)\to T_{p_0}\widetilde M$ by 
$${\widehat f}^F_{p_0,r}(v):=r(F(v)v+({\rm grad}\,F)_v)\,\,\,\,(v\in S^n(1)).$$
Clearly we have $f^F_{p_0,r}={\rm Exp}\circ{\widehat f}^F_{p_0,r}$.  

\vspace{0.5truecm}

\noindent
{\bf Assumption.} 
Let $r$ be such a sufficiently small positive constant as $t^F_r(p_0)$ does not intersect with the cut 
locus of $p_0$.  

\vspace{0.5truecm}

Assume that $\widetilde M$ is a symmetric space of compact type or non-compact type.  
Let $G$ be the identity component of the isometry group of $\widetilde M$ and $G_{p_0}$ 
the isotropy group of $G$ at $p_0$.  For simplicity, set $K:=G_{p_0}$.  
Also, let $\mathfrak g$ and $\mathfrak k$ the Lie algebras of 
$G$ and $K$, respectively, and $\mathfrak g=\mathfrak k+\mathfrak p$ be the canonical decomposition.  
The space $\mathfrak p$ is identified with $T_{p_0}\widetilde M$.  
Fix $v\in\mathfrak p$.  Take a maximal abelian subspace $\mathfrak a_v$ of $\mathfrak p$ 
containing $v$, where ``abelian" means that $R(w_1,w_2)=0$ holds for any elements $w_1$ 
and $w_2$ of $\mathfrak a_v$.  
Then it is shown that $R(w):=R(\cdot,w)w$'s ($w\in\mathfrak a_v$) are simultaneously 
diagonalizable.  
Hence they have common eigenspace decomposition.  
Let $\mathfrak p=\mathfrak a_v\oplus\left(\mathop{\oplus}_{i=1}^k\mathfrak p^v_i\right)$ be their 
common eigenspace decomposition.  
It is clear that there exist an element $\alpha^v_i$ of the dual space 
$\mathfrak a_v^{\ast}$ of $\mathfrak a_v$ such that, for each $w\in\mathfrak a_v$, 
$R(w)=\varepsilon\alpha^v_i(w)^2\,{\rm id}$ holds on $\mathfrak p^v_i$, where 
${\rm id}$ is the identity transformation of $\mathfrak p$ and 
$\varepsilon=1$ (resp. $\varepsilon=-1$) in the case where $G/K$ is of compact type 
(resp. of non-compact type).  
Note that $\alpha^v_i$ is unique up to the $(\pm 1)$-multiple.  
Set $\triangle^v:=\{\pm\alpha^v_1,\cdots,\pm\alpha^v_k\}$.  
For convenience, we denote $\mathfrak p^v_i$ (resp. $\mathfrak a^v$) by 
$\mathfrak p^v_{\alpha_i}$ (resp. $\mathfrak p^v_0$).  
The sysytem $\triangle^v$ gives a root system and 
it is isomorphic to the (restricted) root system of the symmetric pair $(G,K)$.  
Hence, if $\alpha,\beta\in\triangle^v$ and if $\beta=a\alpha$ for some constant $a$, then 
$a=\pm 1\,\,{\rm or}\,\,\pm 2$.  

Denote by $\xi$ and $\nu$ the outward unit normal vector field and the Gauss map of $f^F_{p_0,r}$, 
respectively.  For $\nu$, we have the following fact.  

\vspace{0.5truecm}

\noindent
{\bf Lemma 4.1.} {\sl Assume that $\widetilde M$ is a symmetric space of compact type or 
non-compact type.  For any $v\in S^n(1)$, $\nu(v)=v$ holds.}

\vspace{0.5truecm}

\noindent
{\it Proof.} For simplicity, set $f:=f^F_{p_0,r}$ and $\widehat f:=\widehat f^F_{p_0,r}$.  
Take $X\in T_vS^n(1)$.  Let $c(t)$ ($-\varepsilon<t<\varepsilon$) be 
a curve in $S^n(1)$ with $c'(0)=X$.  
Then we have 
$$\begin{array}{l}
\displaystyle{f_{\ast}X=\left.\frac{d}{dt}\right\vert_{t=0}f(c(t))}\\
\displaystyle{=r{\rm Exp}_{\ast\widehat f(v)}
\left(\left.\frac{d}{dt}\right\vert_{t=0}(F(c(t))c(t)+({\rm grad}\,F)_{c(t)})\right)}\\
\displaystyle{=r{\rm Exp}_{\ast\widehat f(v)}
\left((XF)v+F(v)X+(\nabla^0)^{\iota}_X{\rm grad}\,F\right),}
\end{array}\leqno{(4.1)}$$
where $\nabla^0$ is the Euclidean connection of $T_{p_0}\widetilde M$, 
$\iota$ is the inclusion map of $S^n(1)$ into $T_{p_0}\widetilde M$ and 
$(\nabla^0)^{\iota}$ is the covariant derivative along $\iota$ induced from 
$\nabla^0$.  
On the other hand, since $F$ is $\Phi$-invariant 
(hence invariant with respect to the linear isotropy action 
$K\curvearrowright\mathfrak p$), we may assume that $({\rm grad}\,F)_v\in\mathfrak a_v$ 
by retaking $\mathfrak a_v$ if necessary.  
First we consider the case of $X\in\mathfrak a_v$.  
Then we have $\nabla^0_X{\rm grad}\,F\in\mathfrak a_v$.  
Hence, since $(XF)v+F(v)X+\nabla^0_X{\rm grad}\,F$ belongs to $\mathfrak a_v$, we can derive 
$${\rm Exp}_{\ast\widehat f(v)}((XF)v+F(v)X+\nabla^0_X{\rm grad}\,F)
=\tau_{f(v)}((XF)v+F(v)X+\nabla^0_X{\rm grad}\,F).$$
From $(4.1)$ and this relation, we have 
$$f_{\ast}X=r\tau_{f(v)}((XF)v+F(v)X+\nabla^0_X{\rm grad}\,F).$$
Hence we obtain 
$$\langle f_{\ast}X,\,\tau_{f(v)}(v)\rangle
=r\langle(XF)v+F(v)X+\nabla^0_X{\rm grad}\,F,\,v\rangle=r(XF-\langle{\rm grad}\,F,\,X\rangle)=0.
\leqno{(4.2)}$$
Next we consider the case of $X\in\mathfrak a_v^{\perp}$.  
Take a curve $\hat c:(-\varepsilon,\varepsilon)\to T_{p_0}\widetilde M$ with $\hat c(0)=\widehat f(v)$ 
and ${\hat c}'(0)=(XF)v+F(v)X+\nabla^0_X{\rm grad}\,F$, 
where $\varepsilon$ is a small positive number.  
Define a map $\delta:(-\varepsilon,\varepsilon)\times[0,1]\to\widetilde M$ by 
$$\delta(t,s):={\rm Exp}(s\hat c(t))\,\,
((t,s)\in(-\varepsilon,\varepsilon)\times[0,1]).$$
Set $Y:=\delta_{\ast}\left(\left.\frac{\partial}{\partial t}\right\vert_{t=0}\right)$.  
This vector field $Y$ is the Jacobi field along $\gamma_{f(v)}$ with $Y(0)=0$ and 
$Y'(0)=(XF)v+F(v)X+\nabla^0_X{\rm grad}\,F$.  
Hence it follows from $(2.4)$ that 
$$\begin{array}{l}
\displaystyle{{\rm Exp}_{\ast\widehat f(v)}((XF)v+F(v)X+\nabla^0_X{\rm grad}\,F)}\\
\displaystyle{=Y(1)=\tau_{f(v)}(D^{si}_{\widehat f(v)}((XF)v+F(v)X+\nabla^0_X{\rm grad}\,F)).}
\end{array}$$
On the other hand, from $X\in\mathfrak a_v^{\perp}$, we have $XF=0$ and 
$\nabla^0_X{\rm grad}\,F\in\mathfrak a_v^{\perp}$.  
From these facts, we can show that 
$${\rm Exp}_{\ast\widehat f(v)}((XF)v+F(v)X+\nabla^0_X{\rm grad}\,F)
\in\tau_v(\mathfrak a_v^{\perp}).$$
Hence the relation $(4.2)$ follows from $(4.1)$.  Thus, in both cases, we can derive $(4.2)$.  
Therefore, from the arbitrariness of $X$, we obtain 
$\tau_{f(v)}(v)=\xi_v$, that is, $\nu(v)=v$.
\hspace{0.15truecm}q.e.d.

\vspace{0.5truecm}

Denote by $r_{\widetilde M}$ the first conjugate radius of $\widetilde M$.  
For the anisotropic geodesic sphere, we have the following fact.  

\vspace{0.5truecm}

\noindent
{\bf Theorem 4.2.} {\sl Assume that $\widetilde M$ is an irreducible symmetric space of compact type or 
non-compact type and that it is of rank greater than one.  
Then the following statements ${\rm(i)}\sim{\rm(v)}$ hold:

{\rm(i)} If $r<\frac{r_{\widetilde M}}{2\max_{v\in S^n(1)}\vert\vert F(v)v+({\rm grad}\,F)_v\vert\vert}$, 
then $t_r^F(p_0)$ is anisotropic convex.  

{\rm(ii)} $t_r^F(p_0)$ does not have constant anisotropic principal curvatures.

{\rm(iii)} $t_r^F(p_0)$ is not anisotropic isoparametric.  

{\rm(iv)} If $\widetilde M$ is of compact type, then $t_r^F(p_0)$ is not anisotropic equifocal.  

{\rm(v)} If $\widetilde M$ is of non-compact type, then $t_r^F(p_0)$ is anisotropic equifocal.  
}

\vspace{0.5truecm}

\noindent
{\it Proof.} 
For simplicity, set $f:=f^F_{p_0,r}$ and $\widehat f:=\widehat f^F_{p_0,r}$.  
Denote by $A$ the shape operator of $f$ (for $\xi$).  Also, denote by $\xi_F$ 
the anisotropic transversal vector field of $f$ and $A^F$ the anisotropic shape operator of $f$.  
Take $v\in S^n(1)$.  According to Lemma 4.1, we have $\nu(v)=v$.  
From this fact, it is easy to show that $\gamma'_{f(v)}(1)=r(\xi_F)_v$.  
For simplicity, denote $\gamma_{f(v)}$ by $\gamma$.  
Take $X\in T_v(S^n(1))$ and let $c:(-\varepsilon,\varepsilon)\to S^n(1)$ be 
a curve with $c'(0)=X$, where $\varepsilon$ is a small positive number.  
Define a map $\delta:(-\varepsilon,\varepsilon)\times[0,1]\to\widetilde M$ by 
$$\delta(t,s):={\rm Exp}(s\widehat f(c(t)))\,\,
((t,s)\in(-\varepsilon,\varepsilon)\times[0,1]).$$
Set $Y:=\delta_{\ast}\left(\left.\frac{\partial}{\partial t}\right\vert_{t=0}\right)$, which 
is a Jacobi field along $\gamma$.  We have $Y(0)=0$ and 
$$\begin{array}{l}
\displaystyle{Y'(0)=\left.\frac{d}{dt}\right\vert_{t=0}r\left(
F(c(t))c(t)+({\rm grad}\,F)_{c(t)}\right)}\\
\hspace{1truecm}\displaystyle{=r\left((XF)v+F(v)X+(\nabla^0)^{\iota}_X{\rm grad}\,F\right)}\\
\hspace{1truecm}\displaystyle{=r\left(F(v)X+\nabla^S_X{\rm grad}\,F\right).}
\end{array}$$
So, according to $(2.4)$, we have 
$$Y(s)=\tau_{\gamma\vert_{[0,s]}}\left(srD^{si}_{sr(\xi_F)_v}
(F(v)X+\nabla^S_X{\rm grad}\,F)\right).\leqno{(4.3)}$$
Hence we have 
$$f_{\ast}X=Y(1)=r\tau_{\gamma}\left(D^{si}_{r(\xi_F)_v}
(F(v)X+\nabla^S_X{\rm grad}\,F)\right).\leqno{(4.4)}$$
On the other hand, we have 
$$\begin{array}{l}
\displaystyle{f_{\ast}(A^F_v(X))=-{\widetilde{\nabla}}^f_X\xi_F
=-\frac{1}{r}\left.\left({\widetilde{\nabla}}^{\delta}_{\frac{\partial}{\partial t}}\delta_{\ast}
\left(\frac{\partial}{\partial s}\right)\right)\right\vert_{t=0,s=1}}\\
\hspace{1.95truecm}\displaystyle{=-\frac{1}{r}\left.\left(
{\widetilde{\nabla}}^{\delta}_{\frac{\partial}{\partial s}}
\delta_{\ast}\left(\frac{\partial}{\partial t}\right)\right)\right\vert_{t=0,s=1}
=-\frac{1}{r}Y'(1)}\\
\hspace{1.95truecm}\displaystyle{=-\tau_{\gamma}\left(D^{co}_{r(\xi_F)_v}
(F(v)X+\nabla^S_X{\rm grad}\,F)\right).}
\end{array}\leqno{(4.5)}$$
Let $\mathfrak a_v,\,\triangle^v$ and $\mathfrak p^v_{\alpha}$ be as above.  
Since $F$ is $\Phi$-invariant (hence invariant with respect to the linear isotropy action 
tion $K\curvearrowright\mathfrak p$), we may assume that $({\rm grad}\,F)_v\in\mathfrak a_v$ 
retaking $\mathfrak a_v$ if necessary.  

First we consider the case of $X\in\mathfrak a_v\ominus{\rm Span}\{v\}$.  
Then, since $({\rm grad}\,F)_w\in\mathfrak a_v$ for any $w\in S^n(1)\cap\mathfrak a_v$, 
we have $\nabla^S_X{\rm grad}\,F\in\mathfrak a_v$.  Hence it follows from $(4.4)$ that 
$$f_{\ast}X=r\tau_{\gamma}\left(F(v)X+\nabla^S_X{\rm grad}\,F\right).\leqno{(4.6)}$$
Also it follows from $(4.5)$ that 
$$f_{\ast}(A^F_vX)=-\tau_{\gamma}\left(F(v)X+\nabla^S_X{\rm grad}\,F\right).\leqno{(4.7)}$$
Therefore, we obtain 
$$A^F_vX=-\frac{1}{r}X.\leqno{(4.8)}$$
Also, we have 
$$\left(D^{co}_{s(\xi_F)_v}\circ f_{\ast}
-s(D^{si}_{s(\xi_F)_v}\circ f_{\ast}\circ A^F_v)\right)(X)
=\left(1+\frac{s}{r}\right)f_{\ast}X.\leqno{(4.9)}$$

Next we consider the case of $X\in\mathfrak p^v_{\alpha}$ ($\alpha\in\triangle^v$).  
Then there exist a one-parameter subgroup $\{k_t\}_{t\in{\Bbb R}}$ in $K$ such that 
$\displaystyle{\left.\frac{d}{dt}\right\vert_{t=0}k_t\cdot v=X}$.  
Since $k_{t\ast}\vert_{T_vS^n(1)}$ coincides with the parallel translation along 
$t\mapsto k_t\cdot v$ in $S^n(1)$ and since ${\rm grad}\,F$ is $K$-invariant, 
we obtain $\nabla^S_X{\rm grad}\,F=0$.  Hence it follows from $(4.4)$ that 
$$f_{\ast}X=\frac{F(v)\sin(r\sqrt{\varepsilon}\alpha(\tau_{f(v)}^{-1}((\xi_F)_v)))}
{\sqrt{\varepsilon}\alpha(\tau_{f(v)}^{-1}((\xi_F)_v))}\tau_{\gamma}(X).
\leqno{(4.10)}$$
Also it follows from $(4.5)$ that 
$$f_{\ast}(A^F_vX)=-F(v)\cos(r\sqrt{\varepsilon}\alpha(\tau_{f(v)}^{-1}((\xi_F)_v)))\tau_{\gamma}(X).
\leqno{(4.11)}$$
Therefore, we obtain 
$$A^F_vX=-\frac{\sqrt{\varepsilon}\alpha(\tau_{f(v)}^{-1}((\xi_F)_v))}
{\tan(r\sqrt{\varepsilon}\alpha(\tau_{f(v)}^{-1}((\xi_F)_v)))}X,\leqno{(4.12)}$$
where $\frac{\sqrt{\varepsilon}\alpha(\tau_{f(v)}^{-1}((\xi_F)_v))}
{\tan(r\sqrt{\varepsilon}\alpha(\tau_{f(v)}^{-1}((\xi_F)_v)))}$ means $0$ if 
$\widetilde M$ is of compact type and if 
$r\sqrt{\varepsilon}\alpha(\tau_{f(v)}^{-1}((\xi_F)_v))=\pm\frac{\pi}{2}$.  
Here, when $\widetilde M$ is of compact type, we note that 
$r\sqrt{\varepsilon}\vert\alpha(\tau_{f(v)}^{-1}((\xi_F)_v))\vert$ is smaller than $\pi$ 
because $W(r)$ does not intersect with $C$ (i.e., 
$\vert\vert\widehat f(v)\vert\vert=r\vert\vert(\xi_F)_v\vert\vert$ is smaller than $r_{\widetilde M}$) 
and $r_{\widetilde M}<\frac{\pi\vert\vert(\xi_F)_v\vert\vert}
{\vert\alpha(\tau_{f(v)}^{-1}((\xi_F)_v))\vert}$.  
Also, we have 
$$\begin{array}{l}
\displaystyle{\left(D^{co}_{s(\xi_F)_v}\circ f_{\ast}
-s(D^{si}_{s(\xi_F)_v}\circ f_{\ast}\circ A^F_v)\right)(X)}\\
\displaystyle{=\left(\cos(s\sqrt{\varepsilon}\alpha(\tau_{f(v)}^{-1}((\xi_F)_v)))
+\frac{\sin(s\sqrt{\varepsilon}\alpha(\tau_{f(v)}^{-1}((\xi_F)_v)))}
{\tan(r\sqrt{\varepsilon}\alpha(\tau_{f(v)}^{-1}((\xi_F)_v)))}\right)f_{\ast}X.}
\end{array}\leqno{(4.13)}$$
According to $(4.8)$ and $(4.12)$, the spectrum ${\rm Spec}\,A^F_v$ of $A^F_v$ is given by 
$${\rm Spec}\,A^F_v=\left\{
\left.-\frac{\sqrt{\varepsilon}\alpha(\tau_{f(v)}^{-1}((\xi_F)_v))}
{\tan(r\sqrt{\varepsilon}\alpha(\tau_{f(v)}^{-1}((\xi_F)_v)))}\,\right\vert\,
\alpha\in\triangle^v_+\right\}\cup\left\{-\frac{1}{r}\right\},\leqno{(4.14)}$$
where $\triangle^v_+$ is the positive root system of $\triangle^v$ under some lexicographic ordering of 
$\mathfrak a_v^{\ast}$.  
Since $\widetilde M$ is an irreducible symmetric space of rank greater than one, 
$\alpha((\xi_F)_v)$ depends on the choice of $v\in S^n(1)$.  
Therefore $t_r^F(p_0)$ does not have constant anisotropic principal curvatures and it is not 
anisotropic isoparametric.  Thus the statements (ii) and (iii) follow.  
In particular, if 
$r<\frac{r_{\widetilde M}}{2\max_{v\in S^n(1)}\vert\vert F(v)v+({\rm grad}\,F)_v\vert\vert}$, 
then $r\sqrt{\varepsilon}\vert\alpha(\tau_{f(v)}^{-1}((\xi_F)_v))\vert$ is smaller than $\frac{\pi}{2}$ 
in the case where $\widetilde M$ is of compact type.  
Hence, in this case, it follows from $(4.14)$ that $W_F(r)$ is anisotropic convex.  
Thus the statement (i) follows.  
From $(4.9)$ and $(4.13)$, the set ${\cal AFR}_v$ of all anisotropic focal radii of $t_r^F(p_0)$ at $v$ 
is given by 
$${\cal AFR}_v=
\left\{
\begin{array}{ll}
\displaystyle{\left\{\left.-r+\frac{j\pi}{\alpha(\tau_{f(v)}^{-1}((\xi_F)_v))}
\,\,\right\vert\,\,\alpha\in\triangle_+^v,\,\,j\in{\Bbb Z}\right\}} & 
\displaystyle{(\widetilde M\,:\,{\rm compact}\,\,{\rm type})}\\
\displaystyle{\{-r\}} & \displaystyle{(\widetilde M\,:\,{\rm non-compact}\,\,{\rm type}).}
\end{array}
\right.$$
Hence we obtain the statement (v).  Also, since $\alpha((\xi_F)_v)$ depends on the choice of 
$v\in S^n(1)$, we obtain the statement (iv).  
\hspace{7.75truecm}q.e.d.

\vspace{0.5truecm}

\noindent
{\it Remark 4.1.} If $\widetilde M$ is irreducible and of rank greater than one, then 
geodesic spheres in $\widetilde M$ are not isoparametric and they have not constant principal 
curvatures.  On the basis of this fact, we can conjecture the statements (ii) and (iii) of this theorem 
in advance.  

\section{Anisotropic tubes over certain kind of reflective submanifolds in a symmetric space}
We use the notations in Sections 1-4.  
Let $\widetilde M$ be a symmetric space of compact type or non-compact type.  
In this section, we shall show that anisotropic tubes over a reflective singular orbit of 
a Hermann action of cohomogeneity one on $\widetilde M$ are anisotropic equifocal and 
anisotropic isoparametric hypersurface (for ${\cal F}$).  
Fix $p_0\in\widetilde M$.  Let $G$ and $K$ be as in the previous section.  
Let $\theta$ be the involution of $G$ with 
$({\rm Fix}\,\theta)_0\subset K\subset{\rm Fix}\,\theta$, 
where ${\rm Fix}\,\theta$ is the fixed point group of $\theta$ and 
$({\rm Fix}\,\theta)_0$ is the identity component of ${\rm Fix}\,\theta$.  
Here, in the case where $\widetilde M$ is of compact type (resp. non-compact type), we give 
$\widetilde M$ the $G$-invariant metric induced from the $-\langle\,\,,\,\,\rangle_{\cal K}$ 
(resp. $\langle\,\,,\,\,\rangle_{\cal K}$), where $\langle\,\,,\,\,\rangle_{\cal K}$ is the Killing form 
of the Lie algebra of $G$.  
Let $H$ be a symmetric subgroup of $G$ (i.e., 
$({\rm Fix}\,\tau)_0\subset H\subset{\rm Fix}\,\tau$ for some involution $\tau$ of $G$).  
The natural action of $H$ on $\widetilde M$ is called a {\it Hermann action} 
(see [HPTT], [Kol], [Koi2]).  
In the sequel, we assume that the $H$-action is of cohomogeneity one and commutative, where 
``commutative" means that $\theta\circ\tau=\tau\circ\theta$.  
Let $\mathfrak g,\mathfrak k$ and $\mathfrak h$ be the Lie algebras of $G,K$ 
and $H$, respectively.  We denote the involutions of $\mathfrak g$ induced from 
$\theta$ and $\tau$ by the same symbols $\theta$ and $\tau$, respectively.  
Set $\mathfrak p:={\rm Ker}(\theta+{\rm id})$ and $\mathfrak q:=
{\rm Ker}(\tau+{\rm id})$.  For simplicity, set $B:=Hp_0$, which is reflective because 
$\theta\circ\tau=\tau\circ\theta$.  
Note that $\mathfrak p, \mathfrak p\cap\mathfrak h$ and $\mathfrak p\cap\mathfrak q$ 
are identified with $T_{p_0}\widetilde M,\,T_{p_0}B$ and $T^{\perp}_{p_0}B$, respectively.  
Let $\exp(:T\widetilde M\to\widetilde M)$ be the exponential map of $\widetilde M$ and $\exp^G$ be 
the exponential map of $G$.  
From $\theta\circ\tau=\tau\circ\theta$, it follows that 
$\mathfrak p=\mathfrak p\cap\mathfrak h\oplus \mathfrak p\cap\mathfrak q$.  
We define a map ${\widehat f}^F_{B,r}:T^{\perp_1}B\to T\widetilde M$ by 
$${\widehat f}^F_{B,r}(v):=r\left(\widetilde F(v)v
+({\rm grad}\,\widetilde F\vert_{S^n(1)_{\pi_B(v)}})_v\right)\quad\,\,
(v\in T^{\perp_1}B).$$
Then we have $f^F_{B,r}={\rm exp}\circ{\widehat f}^F_{B,r}$.  
Denote $\nabla^{\perp_B}$ the normal connection of $B$ and 
$P^{hol}_v$ the holonomy subbundle of $T\widetilde M$ through $v$.  
Then we can show the following fact.  

\vspace{0.5truecm}

\noindent
{\bf Lemma 5.1.} {\sl The holonomy invariant elliptic parametric Lagrangian $\widetilde F$ is constant 
over $T^{\perp_1}B$.}

\vspace{0.5truecm}

\noindent
{\it Proof.} 
Take any curve $c$ in $B$.  Since $B$ is totally geodesic, 
the parallel translations along $c$ with respect to $\nabla^{\perp_B}$ and $\widetilde{\nabla}$ 
coincide with each other.  
Hence $P^{hol}_v$ is included by $T^{\perp_1}B$ for each $v\in T^{\perp_1}B$.  
On the other hand, since the $H$-action on $\widetilde M$ is of cohomogeneity one, 
we can derive that the fibre $\pi_B^{-1}(p_0)$ is an orbit of the subaction by $H\cap K$ of 
the linear isotropy group action $K\curvearrowright T_{p_0}\widetilde M$ 
(which is the holonomy group action of $\widetilde M$ at $p_0$).  
Hence we obtain $T^{\perp_1}B=P^{hol}_v$ for any $v\in\pi_B^{-1}(p_0)$.  
Therefore, $\widetilde F$ is constant over $T^{\perp_1}B$ because $\widetilde F$ is holonomy invariant.  
\hspace{10.15truecm}q.e.d.

\vspace{0.5truecm}

Denote by $X_v^L$ the horizontal lift of $X(\in T_p\widetilde M)$ to 
$v(\in T_p\widetilde M)$, where $v$ is regarded as a point of the fibre of $T\widetilde M$ over $p$.  
Also, denote by $c_v^L$ the horizontal lift of a curve $c:[0,\varepsilon)\to\widetilde M$ 
to $v(\in T_{c(0)}\widetilde M)$, where $v$ is regarded as a point of the fibre of $T\widetilde M$ 
over $c(0)$.  Then we can show the following fact.  

\vspace{0.5truecm}

\noindent
{\bf Lemma 5.2.} {\sl Let $X$ be a tangent vector of $B$ at $p$, $\gamma_X$ the geodesic in $B$ with 
$\gamma'_X(0)=X$ and $v$ be an element of $T_p^{\perp_1}B$.  
Then we have 
$$\left.\widetilde{\nabla}^{\gamma_X}_{\frac{d}{dt}}\right\vert_{t=0}
{\widehat f}^F_{B,r}((\gamma_X)_v^L(t))=0.$$}

\vspace{0.5truecm}

\noindent

{\it Proof.} From the holonomy invariantness of $\widetilde F$ and 
$(\gamma_X)_v^L(t)=\tau_{\gamma_X(t)}(v)$, we can derive 
$$\left.\widetilde{\nabla}^{\gamma_X}_{\frac{d}{dt}}\right\vert_{t=0}
(({\rm grad}(\widetilde F\vert_{S^n(1)_{\gamma_X(t)}})_{(\gamma_X)_v^L(t)}))=0.$$
From this relation and Lemma 5.1, the desired relation follows directly.\hspace{1.82truecm}q.e.d.

\vspace{0.5truecm}

Denote by ${\cal H}$ (resp. ${\cal V}$) be the horizontal (resp. vertical) 
distribution on $T^{\perp_1}B$, where ``horizontality" means that so is with respect to 
$\nabla^{\perp_B}$.  
Then we can show the following fact.  

\vspace{0.5truecm}

\noindent
{\bf Theorem 5.3.} {\sl Assume that 
$r<\frac{r_{\widetilde M}}{2\max_{v\in S^n(1)}\vert\vert F(v)v+({\rm grad}\,F)_v\vert\vert}$.  
Then the following statements {\rm(i)}-{\rm(iii)} hold:

{\rm(i)} $f^F_{B,r}$ is an embedding, and $t_r^F(B)$ is an anisotropic equifocal and anisotropic 
isoparametric hypersurface,

{\rm(ii)} $t_r^F(B)$ has constant anisotropic principal curvatures.}

\vspace{0.5truecm}

\noindent
{\it Proof.} For simplicity, set $\widehat f:={\widehat f}^F_{B,r},\,f:=f^F_{B,r}$ 
and $M:=t_r^F(B)$.  Denote by $\xi,\,\xi_F$ and $A^F$ the unit normal vector field , 
the anisotropic normal vector field and the anisotropic shape operator of $f$, respectively.  
For any $p\in B$, denote by $G_p$ the isotropy group of $G$ at $p$.  
For convenience, set $K^p:=G_p$.  
Let $\mathfrak k^p$ be the Lie algebra $K^p$ and $\mathfrak g=\mathfrak k^p+\mathfrak p^p$ be the 
canonical decomposition associated with the symmetric pair $(G,K^p)$.  
Let $\gamma^p_q$ be a shortest geodesic in $\widetilde M$ with $\gamma^p_q(0)=p$ and $\gamma^p_q(1)=q$, 
and $\tau^p_q$ the parallel translation along $\gamma^p_q$.  
Take any $v\in\mathfrak p^p$.  Let $\mathfrak a_v$ be a maximal abelian subspace of $\mathfrak p^p$ 
containing $v$.  Let $\mathfrak p^p=\mathfrak a_v\oplus
\left(\mathop{\oplus}_{i=1}^k\mathfrak p^v_i\right)$ be the common eigenspace decomposition of 
$R(w):=R(\cdot,w)w$'s ($w\in\mathfrak a_v$) and 
$\alpha^v_i$ of the dual space $\mathfrak a_v^{\ast}$ of $\mathfrak a_v$ such that, for each 
$w\in\mathfrak a_v$, $R(w)=\varepsilon\alpha^v_i(w)^2\,{\rm id}$ holds on $\mathfrak p^v_i$, where 
$\varepsilon=1$ (resp. $\varepsilon=-1$) in the case where $\widetilde M$ is of compact type 
(resp. of non-compact type).  
Set $\triangle^v:=\{\pm\alpha^v_1,\cdots,\pm\alpha^v_k\}$.  
For convenience, we denote $\mathfrak p^v_i$ (resp. $\mathfrak a^v$) by 
$\mathfrak p^v_{\alpha_i}$ (resp. $\mathfrak p^v_0$).  

Take any $p\in B$ and any $v\in\pi_B^{-1}(p)$.  
Also, take any $w\in{\cal V}_v$.  
Let $c:(-a,a)\to\pi_B^{-1}(p)$ be a curve with $c'(0)=w$, where $a$ is a positive number.  
Define a map $\delta:(-a,a)\times[0,1]\to\widetilde M$ by 
$\delta(t,s)=\exp_{\pi_B(v)}(s\widehat f(c(t)))\,\,\,((t,s)\in(-a,a)\times[0,1])$.  
Set $\displaystyle{Y:=\left.\frac{\partial\delta}{\partial t}\right\vert_{t=0}}$, 
which is a Jacobi field along $\gamma_{f(v)}$ with $Y(0)=0$.  
Also we have 
$$\begin{array}{l}
\displaystyle{Y'(0)=\left.\frac{d}{dt}\right\vert_{t=0}r\left(
\widetilde F(c(t))c(t)+({\rm grad}(\widetilde F\vert_{S^n(1)_p}))_{c(t)}\right)}\\
\hspace{1truecm}\displaystyle{=r\left((w\widetilde F)v+\widetilde F(v)w
+(\nabla^0)^{\iota}_w{\rm grad}(\widetilde F\vert_{S^n(1)_p})\right)}\\
\hspace{1truecm}\displaystyle{=r\left(\widetilde F(v)w+\nabla^S_w{\rm grad}(\widetilde F\vert_{S^n(1)_p})
\right),}\end{array}$$
where $\nabla^0$ is the Euclidean connection of $T_p\widetilde M$, 
$\iota$ is the inclusion map of $S^n(1)_p$ into $T_p\widetilde M$, 
$(\nabla^0)^{\iota}$ is the covariant derivative along $\iota$ induced from $\nabla^0$ and 
$\nabla^S$ is the Riemannian connection of $S^n(1)_p$.  
Since $\widetilde F$ is holonomy invariant and $w$ belongs to $\mathfrak a_v^{\perp}$, we can derive 
$\nabla^S_w{\rm grad}(\widetilde F\vert_{S^n(1)_p})=0$ (see the proof of Theorem 4.2).  Hence we have 
$Y'(0)=r\widetilde F(v)w$.  Therefore $Y$ is described as 
$$Y(s)=\tau_{\gamma^p_{f(v)}\vert_{[0,s]}}\left(sr\widetilde F(v)D^{si}_{s\widehat f(v)}w\right).
\leqno{(5.1)}$$
Therefore, we obtain 
$$f_{\ast}w=Y(1)=\tau^p_{f(v)}\left(r\widetilde F(v)D^{si}_{\widehat f(v)}w\right).
\leqno{(5.2)}$$

Take $v\in\pi_B^{-1}(p)$ and $X\in T_pB$.  
Let $\gamma_X:(-a,a)\to B$ be the geodesic in $B$ with $\gamma_X'(0)=X$ 
(i.e., $\gamma_X(t)=\exp_p\,tX$) and $\widetilde v$ the parallel unit normal vector field of $B$ 
along $\gamma_X$ with $\widetilde v(0)=v$.  
Define a map $\bar{\delta}:(-a,a)\times[0,1]\to\widetilde M$ by 
$\bar{\delta}(t,s)=\exp_{\gamma_X(t)}(s\widehat f(\widetilde v(t)))\,\,\,
((t,s)\in(-a,a)\times[0,1])$.  
Set $\displaystyle{\bar Y:=\left.\frac{\partial\bar{\delta}}{\partial t}\right\vert_{t=0}}$, 
which is a Jacobi field along $\gamma_{\widehat f(v)}$ with $\bar Y(0)=X$.  
Since $B$ is totally geodesic, we have 
$$\left.{\widetilde{\nabla}}^{\gamma_X}_{\frac{d}{dt}}\right\vert_{t=0}\widetilde v=0,$$
which implies that $\widetilde v(t)=(\gamma_X)_v^L(t)$, that is, 
${\widetilde v}'(0)=X_v^L$.  
On the other hand, since $\widetilde F$ is holonomy invariant, we have $(X^L_v)\widetilde F=0$.  
From these facts and Lemma 5.2, we can derive 
$${\bar Y}'(0)
=\left.{\widetilde{\nabla}}^{\gamma_X}_{\frac{d}{dt}}\right\vert_{t=0}\widehat f(\widetilde v)=0.$$
Therefore $\bar Y$ is described as 
$$\bar Y(s)=\tau_{\gamma^p_{f(v)}\vert_{[0,s]}}
\left(D^{co}_{s\widehat f(v)}X\right).
\leqno{(5.3)}$$
Therefore, we obtain 
$$f_{\ast}(X_v^{L})=\bar Y(1)=\tau^p_{f(v)}(D^{co}_{\widehat f(v)}X).\leqno{(5.4)}$$
Since $r<\frac{r_{\widetilde M}}{2\max_{u\in S^n(1)}\vert\vert F(u)u+({\rm grad}\,F)_u\vert\vert}$ by 
the assumption and $\displaystyle{r_{\widetilde M}<\mathop{\min}_{\alpha\in\triangle^v_+}
\mathop{\min}_{u\in S^n(1)}}\frac{\pi}{\vert\alpha(u)\vert}$, 
we have $\vert\alpha(\widehat f(v))\vert<\frac{\pi}{2}$ for any 
$\alpha\in\triangle^v_+$.  
Therefore, it follows from $(5.2)$ and $(5.4)$ that $f$ is an immersion.  
Furthermore, it is easy to show that $f$ is an embedding.  
From $(5.2)$ and $(5.4)$, we have 
$$\langle f_{\ast}w,\tau^p_{f(v)}(v)\rangle=\langle f_{\ast}(X_v^{L}),\tau^p_{f(v)}(v)\rangle=0.$$
From the arbitrarinesses of $w$ and $X$, these relations imply that $\tau^p_{f(v)}(v)=\xi_v$.  
Furthermore, from this fact, we can derive 
$$(\gamma^p_{f(v)})'(1)=r(\xi_F)_v.\leqno{(5.5)}$$
This fact implies that $f_{-r}(v)=p$.  Hence it follows from the arbitrariness of $v$ that 
$f_{-r}(\pi_B^{-1}(p))=\{p\}$.  Hence we obtain $f_{-r}(T^{\perp_1}B)=B$.  
Thus $B$ is an anisotropic focal submanifold of $f(M)$.  
From $(5.5)$, we have 
$$\begin{array}{l}
\displaystyle{f_{\ast}(A^F_vw)=-\widetilde{\nabla}^f_w\xi_F
=-\frac{1}{r}\left.\left({\widetilde{\nabla}}^{\delta}_{\frac{\partial}{\partial t}}\delta_{\ast}
\left(\frac{\partial}{\partial s}\right)\right)\right\vert_{t=0,s=1}}\\
\hspace{1.6truecm}\displaystyle{=-\frac{1}{r}\left.\left(
{\widetilde{\nabla}}^{\delta}_{\frac{\partial}{\partial s}}\delta_{\ast}\left(\frac{\partial}{\partial t}
\right)\right)\right\vert_{t=0,s=1}
=-\frac{1}{r}Y'(1)}\\
\hspace{1.6truecm}\displaystyle{=-\tau^p_{f(v)}\left(\widetilde F(v)D^{co}_{r(\xi_F)_v}w\right)}
\end{array}\leqno{(5.6)}$$
and 
$$\begin{array}{l}
\displaystyle{f_{\ast}(A^F_v(X_v^{L}))=-\widetilde{\nabla}^f_{X_v^{L}}\xi_F
=-\frac{1}{r}\left.\left({\widetilde{\nabla}}^{\bar{\delta}}_{\frac{\partial}{\partial t}}
\bar{\delta}_{\ast}\left(\frac{\partial}{\partial s}\right)\right)\right\vert_{t=0,s=1}}\\
\hspace{2.17truecm}\displaystyle{=-\frac{1}{r}\left.\left(
{\widetilde{\nabla}}^{\bar{\delta}}_{\frac{\partial}{\partial s}}\bar{\delta}_{\ast}
\left(\frac{\partial}{\partial t}\right)\right)\right\vert_{t=0,s=1}
=-\frac{1}{r}{\bar Y}'(1)}\\
\hspace{2.17truecm}\displaystyle{=\tau^p_{f(v)}\left((({\rm ad}((\xi_F)_v)^2\circ(D^{si}_{r(\xi_F)_v}))X
\right).}
\end{array}\leqno{(5.7)}$$
Assume that $w\in\mathfrak p^v_{\alpha}$.  Then it follows from $(5.2)$ and $(5.6)$ that 
$$A^F_vw=-\frac{\sqrt{\varepsilon}\alpha((\tau^p_{f(v)})^{-1}((\xi_F)_v))}
{\tan(r\sqrt{\varepsilon}\alpha((\tau^p_{f(v)})^{-1}((\xi_F)_v)))}w,\leqno{(5.8)}$$
where $\frac{\sqrt{\varepsilon}\alpha((\tau^p_{f(v)})^{-1}((\xi_F)_v))}
{\tan(r\sqrt{\varepsilon}\alpha((\tau^p_{f(v)})^{-1}((\xi_F)_v)))}$ means $0$ if 
$\widetilde M$ is of compact type and if 
$r\sqrt{\varepsilon}\alpha((\tau^p_{f(v)})^{-1}((\xi_F)_v))=\pm\frac{\pi}{2}$.  
Also, we have 
$$\begin{array}{l}
\displaystyle{\left(D^{co}_{s(\xi_F)_v}\circ f_{\ast}
-s(D^{si}_{s(\xi_F)_v}\circ f_{\ast}\circ A^F_v)\right)(w)}\\
\displaystyle{=\left(\cos(s\sqrt{\varepsilon}\alpha((\tau^p_{f(v)})^{-1}((\xi_F)_v)))
+\frac{\sin(s\sqrt{\varepsilon}\alpha((\tau^p_{f(v)})^{-1}((\xi_F)_v)))}
{\tan(r\sqrt{\varepsilon}\alpha((\tau^p_{f(v)})^{-1}((\xi_F)_v)))}\right)f_{\ast}w.}
\end{array}\leqno{(5.9)}$$
Assume that $X\in\mathfrak p^v_{\alpha}$.  Then it follows from $(5.4)$ and $(5.7)$ that 
$$A^F_v(X_v^{L})=\sqrt{\varepsilon}\alpha((\tau^p_{f(v)})^{-1}((\xi_F)_v))
\tan(r\sqrt{\varepsilon}\alpha((\tau^p_{f(v)})^{-1}((\xi_F)_v))))X_v^{L}.\leqno{(5.10)}$$
Also, we have 
$$\begin{array}{l}
\displaystyle{\left(D^{co}_{s(\xi_F)_v}\circ f_{\ast}
-s(D^{si}_{s(\xi_F)_v}\circ f_{\ast}\circ A^F_v)\right)(X_v^{L})}\\
\displaystyle{=\left(\cos(s\sqrt{\varepsilon}\alpha((\tau^p_{f(v)})^{-1}((\xi_F)_v)))\right.}\\
\hspace{0.8truecm}\displaystyle{\left.-\sin(s\sqrt{\varepsilon}\alpha((\tau^p_{f(v)})^{-1}((\xi_F)_v)))
\tan(r\sqrt{\varepsilon}\alpha((\tau^p_{f(v)})^{-1}((\xi_F)_v)))\right)f_{\ast}(X_v^{L}).}
\end{array}\leqno{(5.11)}$$
According to $(5.8)$ and $(5.10)$, we obtain 
$$\begin{array}{l}
\displaystyle{{\rm Spec}\,A^F_v=\left\{
\left.-\frac{\sqrt{\varepsilon}\alpha((\tau^p_{f(v)})^{-1}((\xi_F)_v))}
{\tan(r\sqrt{\varepsilon}\alpha((\tau^p_{f(v)})^{-1}((\xi_F)_v)))}\,\right\vert\,
\alpha\in\triangle^v_+\,\,{\rm s.t.}\,\,\mathfrak p^v_{\alpha}\cap\mathfrak q\not=\{0\}\right\}}\\
\hspace{2.3truecm}\displaystyle{\cup\left\{
\sqrt{\varepsilon}\alpha((\tau^p_{f(v)})^{-1}((\xi_F)_v))
\tan(r\sqrt{\varepsilon}\alpha((\tau^p_{f(v)})^{-1}((\xi_F)_v)))\right.}\\
\hspace{4.4truecm}\displaystyle{\left.\vert\,\alpha\in\triangle^v_+\cup\{0\}\,\,{\rm s.t.}\,\,
\mathfrak p^v_{\alpha}\cap\mathfrak h\not=\{0\}\right\},}
\end{array}\leqno{(5.12)}$$
where $\mathfrak p^v_0$ means $\mathfrak a_v$.  
Also, according to $(5.9)$ and $(5.11)$, we obtain 
$${\cal AFR}_v=
\left\{
\begin{array}{l}
\displaystyle{\left\{\left.-r+\frac{j\pi}{\alpha((\tau^p_{f(v)})^{-1}((\xi_F)_v))}
\,\,\right\vert\,\,j\in{\Bbb Z},\,\,\,\alpha\in\triangle_+^v\,\,{\rm s.t.}\,\,
\mathfrak p^v_{\alpha}\cap\mathfrak q\not=\{0\}\right\}}\\
\displaystyle{\bigcup\left\{\left.r+\frac{(2j+1)\pi}{2\alpha((\tau^p_{f(v)})^{-1}((\xi_F)_v))}
\,\,\right\vert\,\,j\in{\Bbb Z},\,\,\,\alpha\in\triangle_+^v\,\,{\rm s.t.}\,\,
\mathfrak p^v_{\alpha}\cap\mathfrak h\not=\{0\}\right\}}\\
\hspace{6.3truecm}({\rm when}\,\,\widetilde M\,:\,{\rm compact}\,\,{\rm type})\\
\displaystyle{\{-r\}}\hspace{4truecm}({\rm when}\,\,\widetilde M\,:\,{\rm non-compact}\,\,{\rm type}).
\end{array}
\right.\leqno{(5.13)}$$

Take a geodesic $\bar{\gamma}$ in $B$ with $\bar{\gamma}(0)=p$.  
Set $\bar X:={\bar{\gamma}}'(0)$.  
We have $\bar{\gamma}(t)={\rm Exp}\,t\bar X$ and $\tau_{\bar{\gamma}\vert_{[0,t]}}
=(\exp^G\,t\bar X)_{\ast p}$.  Set $g_t:=\exp^G\,t\bar X$.  
According to Lemma 5.2, we have 
$$\begin{array}{l}
\displaystyle{(\tau^{\bar{\gamma}(t)}_{f(\bar{\gamma}_v^L(t))})^{-1}
((\xi_F)_{\bar{\gamma}_v^L(t)})=\frac{1}{r}\widehat f({\bar{\gamma}}_v^L(t))}\\
\displaystyle{=\frac{1}{r}(g_t)_{\ast p}((\tau^p_{f(\bar{\gamma}_v^L(t))})^{-1}\widehat f(v))}\\
\displaystyle{=(g_t)_{\ast p}((\tau^p_{f(v)})^{-1}((\xi_F)_v)),}
\end{array}$$
where $\tau^{\bar{\gamma}(t)}_{f(\bar{\gamma}_v^L(t))}$ is defined in similar to $\tau^p_q$.  
From this fact, we can derive ${\rm Spec}\,A^F_{\bar{\gamma}_v^L(t)}={\rm Spec}\,A^F_v$ and 
${\cal AFR}_{\bar{\gamma}_v^L(t)}={\cal AFR}_v$.  
Therefore, it follows from the arbitrariness of $\bar{\gamma}$ and $P^{hol}_v=T^{\perp_1}B$ that 
${\rm Spec}\,A^F_{\bullet}$ and ${\cal AFR}_{\bullet}$ are independent of the choice of 
$\bullet\in T^{\perp_1}B$.  
Therefore, $M$ is anisotropic equifocal and has constant anisotropic principal curvatures.  
For each $t$ close to $0$, the anisotropic parallel hypersurface 
$f_t:T^{\perp_1}B\to\widetilde M$ is given by 
$$f_t(v):=\exp_{\pi_B(v)}\left(\left(\frac{t}{r}+1\right)\widehat f(v)\right)\quad\,\,
(v\in T^{\perp_1}B).$$
Set $M_t:=f_t(T^{\perp_1}B)$.  
As above, it is shown that $M_t$'s have constant anisotropic principal curvatures 
and hence they have constant anisotropic mean curvature.  That is, $M$ is anisotropic 
isoparametric.  This completes the proof.  \hspace{5.25truecm}q.e.d.

\vspace{1.4truecm}

\centerline{
\unitlength 0.1in
\begin{picture}( 38.9000, 15.3000)(  0.6000,-23.3000)
%
\special{pn 8}%
\special{pa 1350 1970}%
\special{pa 3150 1970}%
\special{fp}%
%
\special{pn 8}%
\special{pa 3150 1970}%
\special{pa 3950 1570}%
\special{fp}%
%
\special{pn 8}%
\special{pa 1350 1970}%
\special{pa 2150 1570}%
\special{fp}%
%
\special{pn 8}%
\special{ar 2560 1750 420 120  6.2831853 6.2831853}%
\special{ar 2560 1750 420 120  0.0000000 3.1415927}%
%
\special{pn 8}%
\special{ar 2550 1770 424 424  2.7386811 6.2831853}%
\special{ar 2550 1770 424 424  0.0000000 0.3805064}%
%
\special{pn 8}%
\special{ar 2550 1770 416 416  0.5840811 2.5535901}%
%
\special{pn 8}%
\special{pa 2550 1060}%
\special{pa 2550 1850}%
\special{fp}%
%
\special{pn 8}%
\special{pa 2550 1920}%
\special{pa 2550 1940}%
\special{fp}%
%
\special{pn 8}%
\special{pa 2550 2010}%
\special{pa 2550 2330}%
\special{fp}%
%
\special{pn 8}%
\special{pa 3940 1570}%
\special{pa 2960 1570}%
\special{fp}%
%
\special{pn 8}%
\special{pa 2870 1570}%
\special{pa 2580 1570}%
\special{fp}%
%
\special{pn 8}%
\special{pa 2500 1570}%
\special{pa 2230 1570}%
\special{fp}%
%
\special{pn 8}%
\special{pa 2550 1750}%
\special{pa 2970 1750}%
\special{fp}%
\special{sh 1}%
\special{pa 2970 1750}%
\special{pa 2904 1730}%
\special{pa 2918 1750}%
\special{pa 2904 1770}%
\special{pa 2970 1750}%
\special{fp}%
%
\special{pn 8}%
\special{pa 2550 1750}%
\special{pa 3140 1750}%
\special{fp}%
\special{sh 1}%
\special{pa 3140 1750}%
\special{pa 3074 1730}%
\special{pa 3088 1750}%
\special{pa 3074 1770}%
\special{pa 3140 1750}%
\special{fp}%
%
\special{pn 8}%
\special{pa 2550 1750}%
\special{pa 2550 1460}%
\special{fp}%
\special{sh 1}%
\special{pa 2550 1460}%
\special{pa 2530 1528}%
\special{pa 2550 1514}%
\special{pa 2570 1528}%
\special{pa 2550 1460}%
\special{fp}%
%
\special{pn 8}%
\special{pa 2550 1750}%
\special{pa 3130 1460}%
\special{fp}%
\special{sh 1}%
\special{pa 3130 1460}%
\special{pa 3062 1472}%
\special{pa 3082 1484}%
\special{pa 3080 1508}%
\special{pa 3130 1460}%
\special{fp}%
%
\special{pn 8}%
\special{pa 3140 1750}%
\special{pa 3140 1460}%
\special{dt 0.045}%
%
\special{pn 8}%
\special{pa 3140 1460}%
\special{pa 2550 1460}%
\special{dt 0.045}%
%
\special{pn 8}%
\special{ar 2550 1750 420 110  4.7995042 6.2831853}%
%
\special{pn 8}%
\special{ar 2550 1750 410 110  3.1415927 4.6613478}%
%
\special{pn 8}%
\special{pa 1720 1510}%
\special{pa 1940 1850}%
\special{dt 0.045}%
\special{sh 1}%
\special{pa 1940 1850}%
\special{pa 1922 1784}%
\special{pa 1912 1806}%
\special{pa 1888 1806}%
\special{pa 1940 1850}%
\special{fp}%
%
\special{pn 8}%
\special{pa 2370 960}%
\special{pa 2550 1100}%
\special{dt 0.045}%
\special{sh 1}%
\special{pa 2550 1100}%
\special{pa 2510 1044}%
\special{pa 2508 1068}%
\special{pa 2486 1076}%
\special{pa 2550 1100}%
\special{fp}%
%
\special{pn 8}%
\special{pa 3010 2070}%
\special{pa 2800 1750}%
\special{dt 0.045}%
\special{sh 1}%
\special{pa 2800 1750}%
\special{pa 2820 1818}%
\special{pa 2830 1796}%
\special{pa 2854 1796}%
\special{pa 2800 1750}%
\special{fp}%
%
\special{pn 8}%
\special{pa 1960 2110}%
\special{pa 2400 1870}%
\special{dt 0.045}%
\special{sh 1}%
\special{pa 2400 1870}%
\special{pa 2332 1884}%
\special{pa 2354 1896}%
\special{pa 2352 1920}%
\special{pa 2400 1870}%
\special{fp}%
%
\special{pn 8}%
\special{pa 2820 2270}%
\special{pa 2710 2040}%
\special{dt 0.045}%
\special{sh 1}%
\special{pa 2710 2040}%
\special{pa 2722 2110}%
\special{pa 2734 2088}%
\special{pa 2758 2092}%
\special{pa 2710 2040}%
\special{fp}%
\put(16.8000,-15.8000){\makebox(0,0)[rb]{$T^{\perp}_{p_0}B$}}%
\put(23.6000,-9.7000){\makebox(0,0)[rb]{$T_{p_0}B$}}%
\put(26.9000,-11.8000){\makebox(0,0)[lb]{$r({\rm grad}\widetilde F\vert_{T^1_{p_0}\widetilde M})_v\,\,(\in\mathfrak p^v_0\cap\mathfrak h\subset T_{p_0}B)$}}%
\put(19.3000,-21.4000){\makebox(0,0)[rt]{$\pi_B^{-1}(p_0)$}}%
\put(28.0000,-23.1000){\makebox(0,0)[lt]{$S^n(1)$}}%
\put(29.9000,-21.2000){\makebox(0,0)[lt]{$v$}}%
%
\special{pn 8}%
\special{pa 3340 2050}%
\special{pa 3030 1750}%
\special{dt 0.045}%
\special{sh 1}%
\special{pa 3030 1750}%
\special{pa 3064 1812}%
\special{pa 3068 1788}%
\special{pa 3092 1782}%
\special{pa 3030 1750}%
\special{fp}%
\put(33.9000,-20.7000){\makebox(0,0)[lt]{$r\widetilde F(v)v$}}%
\put(31.9000,-14.6000){\makebox(0,0)[lb]{${\widehat f}^F_{B,r}(v)$}}%
%
\special{pn 8}%
\special{pa 3010 1200}%
\special{pa 2550 1610}%
\special{dt 0.045}%
\special{sh 1}%
\special{pa 2550 1610}%
\special{pa 2614 1582}%
\special{pa 2590 1576}%
\special{pa 2586 1552}%
\special{pa 2550 1610}%
\special{fp}%
\end{picture}%
\hspace{3truecm}}

\vspace{1truecm}

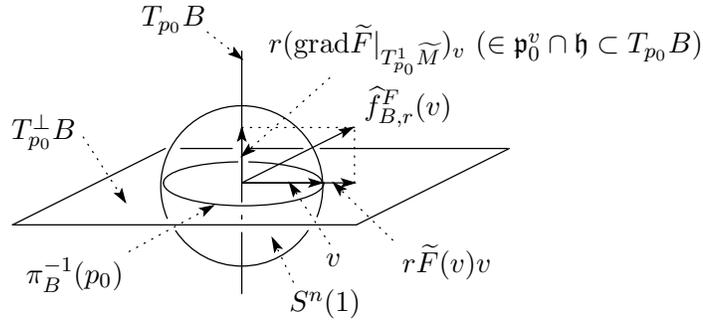
\centerline{{\bf Figure 3.}}

\vspace{1truecm}

\centerline{
\unitlength 0.1in
\begin{picture}( 52.6100, 22.0300)(-13.4600,-24.7000)
%
\special{pn 8}%
\special{ar 1928 2434 1488 314  4.3042258 5.7757269}%
%
\special{pn 8}%
\special{ar 2524 1770 384 120  6.2831853 6.2831853}%
\special{ar 2524 1770 384 120  0.0000000 3.1415927}%
%
\special{pn 8}%
\special{ar 2524 1770 384 130  4.8396967 6.2831853}%
%
\special{pn 8}%
\special{ar 2524 1770 384 130  3.1415927 4.5850812}%
%
\special{pn 8}%
\special{ar 3282 1780 376 1236  2.8263112 3.8891710}%
%
\special{pn 8}%
\special{ar 3224 1852 312 1028  2.4965248 2.7469318}%
%
\special{pn 8}%
\special{ar 2514 1696 376 1236  2.8263112 3.8891710}%
%
\special{pn 8}%
\special{ar 2504 1660 374 1218  2.4974388 2.7469103}%
%
\special{pn 8}%
\special{ar 2092 2392 764 946  3.3837530 4.4468074}%
%
\special{pn 8}%
\special{ar 2236 2418 3178 940  4.6113613 4.6380403}%
%
\special{pn 8}%
\special{ar 2918 1734 374 1236  2.8276677 3.8880495}%
%
\special{pn 8}%
\special{ar 2926 1724 374 1218  2.4964576 2.7469103}%
%
\special{pn 8}%
\special{pa 2944 1238}%
\special{pa 2928 1266}%
\special{pa 2904 1286}%
\special{pa 2876 1300}%
\special{pa 2846 1310}%
\special{pa 2814 1316}%
\special{pa 2782 1322}%
\special{pa 2750 1326}%
\special{pa 2718 1330}%
\special{pa 2686 1332}%
\special{pa 2654 1332}%
\special{pa 2622 1332}%
\special{pa 2590 1332}%
\special{pa 2560 1330}%
\special{pa 2528 1326}%
\special{pa 2496 1324}%
\special{pa 2464 1318}%
\special{pa 2432 1314}%
\special{pa 2400 1308}%
\special{pa 2370 1300}%
\special{pa 2338 1292}%
\special{pa 2308 1282}%
\special{pa 2278 1272}%
\special{pa 2250 1258}%
\special{pa 2222 1242}%
\special{pa 2196 1222}%
\special{pa 2180 1196}%
\special{pa 2178 1180}%
\special{sp}%
%
\special{pn 8}%
\special{pa 2620 1088}%
\special{pa 2652 1092}%
\special{pa 2682 1098}%
\special{pa 2714 1102}%
\special{pa 2746 1108}%
\special{pa 2776 1116}%
\special{pa 2808 1126}%
\special{pa 2838 1136}%
\special{pa 2866 1148}%
\special{pa 2894 1166}%
\special{pa 2918 1184}%
\special{pa 2940 1208}%
\special{pa 2946 1236}%
\special{sp}%
%
\special{pn 8}%
\special{pa 2178 1190}%
\special{pa 2190 1160}%
\special{pa 2214 1140}%
\special{pa 2242 1126}%
\special{pa 2270 1114}%
\special{pa 2302 1106}%
\special{pa 2332 1098}%
\special{pa 2364 1092}%
\special{pa 2396 1090}%
\special{pa 2428 1088}%
\special{pa 2460 1086}%
\special{pa 2492 1086}%
\special{pa 2522 1086}%
\special{sp}%
%
\special{pn 20}%
\special{sh 1}%
\special{ar 2542 1770 10 10 0  6.28318530717959E+0000}%
\special{sh 1}%
\special{ar 2542 1770 10 10 0  6.28318530717959E+0000}%
%
\special{pn 8}%
\special{ar 2544 1580 626 146  6.2831853 6.2831853}%
\special{ar 2544 1580 626 146  0.0000000 3.1415927}%
%
\special{pn 8}%
\special{ar 2544 1580 626 158  4.8397743 6.2831853}%
%
\special{pn 8}%
\special{ar 2544 1580 624 158  3.1415927 4.5853358}%
%
\special{pn 8}%
\special{pa 3176 2158}%
\special{pa 3174 2126}%
\special{pa 3170 2094}%
\special{pa 3166 2062}%
\special{pa 3164 2030}%
\special{pa 3162 1998}%
\special{pa 3158 1966}%
\special{pa 3156 1934}%
\special{pa 3156 1902}%
\special{pa 3154 1870}%
\special{pa 3154 1838}%
\special{pa 3152 1806}%
\special{pa 3152 1774}%
\special{pa 3152 1742}%
\special{pa 3152 1710}%
\special{pa 3152 1678}%
\special{pa 3154 1646}%
\special{pa 3154 1614}%
\special{pa 3156 1582}%
\special{pa 3156 1550}%
\special{pa 3158 1518}%
\special{pa 3160 1486}%
\special{pa 3164 1454}%
\special{pa 3164 1422}%
\special{pa 3168 1392}%
\special{pa 3172 1360}%
\special{pa 3174 1328}%
\special{pa 3178 1296}%
\special{pa 3184 1264}%
\special{pa 3188 1232}%
\special{pa 3192 1200}%
\special{pa 3198 1170}%
\special{pa 3202 1138}%
\special{pa 3208 1106}%
\special{pa 3214 1074}%
\special{pa 3220 1044}%
\special{pa 3226 1012}%
\special{pa 3232 980}%
\special{pa 3240 950}%
\special{pa 3248 918}%
\special{pa 3256 888}%
\special{pa 3266 856}%
\special{pa 3274 826}%
\special{pa 3282 796}%
\special{pa 3294 766}%
\special{pa 3304 734}%
\special{pa 3316 704}%
\special{pa 3326 674}%
\special{pa 3326 672}%
\special{sp}%
%
\special{pn 8}%
\special{pa 3224 1032}%
\special{pa 3210 1060}%
\special{pa 3184 1080}%
\special{pa 3156 1094}%
\special{pa 3126 1106}%
\special{pa 3094 1112}%
\special{pa 3064 1120}%
\special{pa 3032 1126}%
\special{pa 3000 1132}%
\special{pa 2968 1136}%
\special{pa 2936 1138}%
\special{pa 2904 1142}%
\special{pa 2872 1144}%
\special{pa 2840 1146}%
\special{pa 2808 1146}%
\special{pa 2776 1146}%
\special{pa 2744 1146}%
\special{pa 2712 1146}%
\special{pa 2680 1146}%
\special{pa 2648 1146}%
\special{pa 2616 1146}%
\special{pa 2584 1144}%
\special{pa 2554 1142}%
\special{pa 2522 1140}%
\special{pa 2490 1136}%
\special{pa 2458 1134}%
\special{pa 2426 1130}%
\special{pa 2394 1126}%
\special{pa 2362 1122}%
\special{pa 2330 1116}%
\special{pa 2300 1110}%
\special{pa 2268 1106}%
\special{pa 2236 1098}%
\special{pa 2204 1092}%
\special{pa 2174 1084}%
\special{pa 2144 1074}%
\special{pa 2112 1064}%
\special{pa 2084 1052}%
\special{pa 2054 1040}%
\special{pa 2026 1026}%
\special{pa 2000 1006}%
\special{pa 1982 980}%
\special{pa 1978 962}%
\special{sp}%
%
\special{pn 8}%
\special{pa 2686 866}%
\special{pa 2718 870}%
\special{pa 2750 872}%
\special{pa 2782 876}%
\special{pa 2814 882}%
\special{pa 2846 886}%
\special{pa 2876 890}%
\special{pa 2908 896}%
\special{pa 2940 900}%
\special{pa 2970 910}%
\special{pa 3002 916}%
\special{pa 3032 926}%
\special{pa 3062 934}%
\special{pa 3094 944}%
\special{pa 3122 958}%
\special{pa 3152 970}%
\special{pa 3178 988}%
\special{pa 3202 1008}%
\special{pa 3216 1036}%
\special{pa 3216 1046}%
\special{sp}%
%
\special{pn 8}%
\special{pa 1966 968}%
\special{pa 1980 940}%
\special{pa 2004 920}%
\special{pa 2032 906}%
\special{pa 2062 892}%
\special{pa 2092 884}%
\special{pa 2122 876}%
\special{pa 2154 870}%
\special{pa 2186 864}%
\special{pa 2216 858}%
\special{pa 2248 854}%
\special{pa 2280 852}%
\special{pa 2312 848}%
\special{pa 2344 846}%
\special{pa 2376 844}%
\special{pa 2408 844}%
\special{pa 2440 842}%
\special{pa 2472 842}%
\special{pa 2504 842}%
\special{pa 2526 842}%
\special{sp}%
%
\special{pn 8}%
\special{pa 1918 2092}%
\special{pa 1916 2060}%
\special{pa 1914 2028}%
\special{pa 1912 1996}%
\special{pa 1912 1964}%
\special{pa 1910 1932}%
\special{pa 1908 1900}%
\special{pa 1908 1868}%
\special{pa 1908 1836}%
\special{pa 1906 1804}%
\special{pa 1906 1772}%
\special{pa 1906 1740}%
\special{pa 1906 1708}%
\special{pa 1906 1676}%
\special{pa 1906 1644}%
\special{pa 1908 1612}%
\special{pa 1908 1580}%
\special{pa 1910 1548}%
\special{pa 1910 1516}%
\special{pa 1912 1486}%
\special{pa 1914 1454}%
\special{pa 1916 1422}%
\special{pa 1918 1390}%
\special{pa 1920 1358}%
\special{pa 1922 1326}%
\special{pa 1924 1294}%
\special{pa 1928 1262}%
\special{pa 1932 1230}%
\special{pa 1936 1198}%
\special{pa 1938 1166}%
\special{pa 1942 1134}%
\special{pa 1946 1102}%
\special{pa 1950 1070}%
\special{pa 1956 1040}%
\special{pa 1962 1008}%
\special{pa 1966 976}%
\special{pa 1970 944}%
\special{pa 1978 914}%
\special{pa 1982 882}%
\special{pa 1988 850}%
\special{pa 1994 818}%
\special{pa 2002 788}%
\special{pa 2008 756}%
\special{pa 2016 726}%
\special{pa 2024 694}%
\special{pa 2032 664}%
\special{pa 2040 632}%
\special{pa 2050 602}%
\special{pa 2058 572}%
\special{pa 2060 566}%
\special{sp}%
%
\special{pn 8}%
\special{ar 2486 1614 586 1900  2.7525983 2.8485271}%
%
\special{pn 8}%
\special{ar 3628 1576 470 1972  2.6916899 2.7470255}%
%
\special{pn 8}%
\special{pa 3484 842}%
\special{pa 2870 1368}%
\special{dt 0.045}%
\special{sh 1}%
\special{pa 2870 1368}%
\special{pa 2934 1340}%
\special{pa 2910 1334}%
\special{pa 2908 1310}%
\special{pa 2870 1368}%
\special{fp}%
%
\special{pn 8}%
\special{pa 3560 1120}%
\special{pa 3090 1332}%
\special{dt 0.045}%
\special{sh 1}%
\special{pa 3090 1332}%
\special{pa 3158 1322}%
\special{pa 3138 1310}%
\special{pa 3142 1286}%
\special{pa 3090 1332}%
\special{fp}%
%
\special{pn 8}%
\special{pa 1352 1442}%
\special{pa 1678 1922}%
\special{dt 0.045}%
\special{sh 1}%
\special{pa 1678 1922}%
\special{pa 1658 1856}%
\special{pa 1648 1878}%
\special{pa 1624 1878}%
\special{pa 1678 1922}%
\special{fp}%
\put(36.0800,-11.5600){\makebox(0,0)[lb]{$t^F_r(B)$ (outside tube)}}%
\put(35.0200,-7.7800){\makebox(0,0)[lb]{$H{\rm Exp}(rv)$ (inside tube)}}%
%
\special{pn 8}%
\special{ar 3282 1756 740 1688  3.6622686 3.8729035}%
%
\special{pn 8}%
\special{ar 2878 1718 740 1688  3.6617577 3.8729035}%
%
\special{pn 8}%
\special{ar 3646 1810 740 1688  3.6617577 3.8729035}%
%
\special{pn 20}%
\special{sh 1}%
\special{ar 2908 1774 10 10 0  6.28318530717959E+0000}%
\special{sh 1}%
\special{ar 2908 1774 10 10 0  6.28318530717959E+0000}%
%
\special{pn 20}%
\special{sh 1}%
\special{ar 3166 1590 10 10 0  6.28318530717959E+0000}%
\special{sh 1}%
\special{ar 3166 1590 10 10 0  6.28318530717959E+0000}%
\put(15.3400,-13.9500){\makebox(0,0)[rb]{$\exp^{\perp}(T^{\perp}_{p_0}B)$}}%
%
\special{pn 8}%
\special{pa 3512 1368}%
\special{pa 3176 1590}%
\special{dt 0.045}%
\special{sh 1}%
\special{pa 3176 1590}%
\special{pa 3244 1570}%
\special{pa 3222 1560}%
\special{pa 3222 1536}%
\special{pa 3176 1590}%
\special{fp}%
%
\special{pn 8}%
\special{pa 3828 1718}%
\special{pa 2926 1774}%
\special{dt 0.045}%
\special{sh 1}%
\special{pa 2926 1774}%
\special{pa 2994 1790}%
\special{pa 2980 1770}%
\special{pa 2992 1750}%
\special{pa 2926 1774}%
\special{fp}%
\put(39.1500,-16.8100){\makebox(0,0)[lt]{${\rm Exp}(rv)$}}%
\put(35.7900,-13.8600){\makebox(0,0)[lb]{$f^F_{B,r}(v)$}}%
%
\special{pn 8}%
\special{pa 2494 484}%
\special{pa 2696 714}%
\special{dt 0.045}%
\special{sh 1}%
\special{pa 2696 714}%
\special{pa 2668 650}%
\special{pa 2662 674}%
\special{pa 2638 676}%
\special{pa 2696 714}%
\special{fp}%
\put(25.1300,-4.3700){\makebox(0,0)[rb]{$B$}}%
\put(24.7500,-17.6400){\makebox(0,0)[rt]{$p_0$}}%
%
\special{pn 8}%
\special{pa 3210 2294}%
\special{pa 3218 2264}%
\special{pa 3226 2232}%
\special{pa 3238 2202}%
\special{pa 3248 2172}%
\special{pa 3258 2142}%
\special{pa 3272 2112}%
\special{pa 3284 2084}%
\special{pa 3298 2054}%
\special{pa 3314 2028}%
\special{pa 3330 2000}%
\special{pa 3346 1972}%
\special{pa 3364 1944}%
\special{pa 3382 1920}%
\special{pa 3402 1894}%
\special{pa 3422 1868}%
\special{pa 3444 1846}%
\special{pa 3466 1822}%
\special{pa 3488 1798}%
\special{pa 3512 1778}%
\special{pa 3536 1756}%
\special{pa 3560 1736}%
\special{pa 3586 1718}%
\special{pa 3612 1700}%
\special{pa 3640 1682}%
\special{pa 3668 1666}%
\special{pa 3696 1652}%
\special{pa 3724 1638}%
\special{pa 3754 1626}%
\special{pa 3784 1614}%
\special{pa 3814 1604}%
\special{pa 3826 1602}%
\special{sp}%
%
\special{pn 8}%
\special{ar 1948 2650 4464 1162  4.9962037 5.1476051}%
%
\special{pn 8}%
\special{ar 1948 2640 4838 1154  4.9220892 4.9489993}%
%
\special{pn 8}%
\special{ar 1948 2640 4570 1154  4.8580607 4.9170902}%
%
\special{pn 8}%
\special{ar 1966 2704 4924 1218  4.7554505 4.8194527}%
%
\special{pn 8}%
\special{ar 1814 2188 3168 710  4.8002859 4.8149948}%
\end{picture}%
\hspace{6truecm}}

\vspace{0.75truecm}

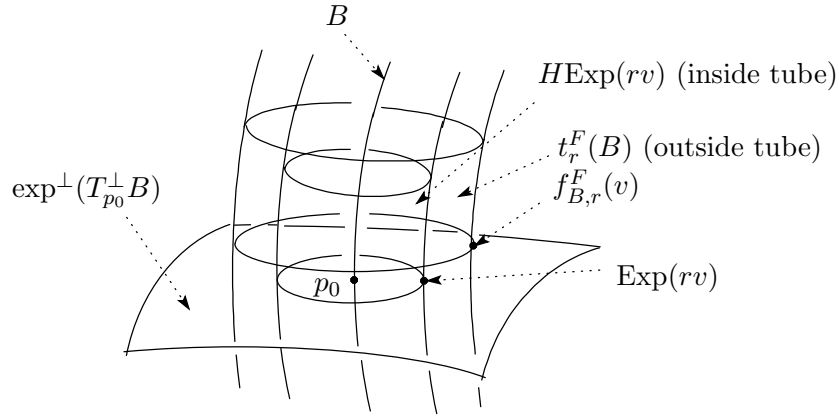
\centerline{{\bf Figure 4.}}

\vspace{1truecm}

\centerline{
\unitlength 0.1in
\begin{picture}( 64.5300, 24.9100)(-15.7000,-27.3700)
%
\special{pn 8}%
\special{ar 1562 2170 1254 276  4.3044905 5.7737703}%
%
\special{pn 8}%
\special{pa 2638 2038}%
\special{pa 2648 2008}%
\special{pa 2660 1978}%
\special{pa 2672 1948}%
\special{pa 2688 1920}%
\special{pa 2702 1892}%
\special{pa 2718 1864}%
\special{pa 2734 1836}%
\special{pa 2752 1810}%
\special{pa 2770 1784}%
\special{pa 2790 1758}%
\special{pa 2810 1732}%
\special{pa 2830 1710}%
\special{pa 2852 1686}%
\special{pa 2876 1664}%
\special{pa 2898 1640}%
\special{pa 2922 1620}%
\special{pa 2948 1600}%
\special{pa 2974 1582}%
\special{pa 3000 1562}%
\special{pa 3026 1544}%
\special{pa 3054 1530}%
\special{pa 3082 1516}%
\special{pa 3112 1500}%
\special{pa 3142 1490}%
\special{pa 3170 1478}%
\special{pa 3202 1468}%
\special{pa 3232 1460}%
\special{pa 3262 1450}%
\special{pa 3294 1446}%
\special{sp}%
%
\special{pn 8}%
\special{pa 1078 1930}%
\special{pa 1082 1898}%
\special{pa 1088 1866}%
\special{pa 1094 1836}%
\special{pa 1102 1804}%
\special{pa 1110 1772}%
\special{pa 1120 1742}%
\special{pa 1130 1712}%
\special{pa 1144 1684}%
\special{pa 1156 1654}%
\special{pa 1170 1626}%
\special{pa 1184 1596}%
\special{pa 1202 1570}%
\special{pa 1218 1542}%
\special{pa 1236 1516}%
\special{pa 1258 1492}%
\special{pa 1278 1468}%
\special{pa 1300 1444}%
\special{pa 1322 1422}%
\special{pa 1346 1400}%
\special{pa 1372 1380}%
\special{pa 1398 1362}%
\special{pa 1426 1346}%
\special{pa 1452 1330}%
\special{pa 1474 1322}%
\special{sp}%
%
\special{pn 8}%
\special{ar 2226 2018 2046 680  4.9213420 5.2614898}%
%
\special{pn 8}%
\special{ar 2394 1554 316 1084  2.8259026 3.8898760}%
%
\special{pn 20}%
\special{sh 1}%
\special{ar 2080 1586 10 10 0  6.28318530717959E+0000}%
\special{sh 1}%
\special{ar 2080 1586 10 10 0  6.28318530717959E+0000}%
%
\special{pn 8}%
\special{ar 2072 1388 526 128  6.2831853 6.2831853}%
\special{ar 2072 1388 526 128  0.0000000 3.1415927}%
%
\special{pn 8}%
\special{ar 2064 1396 528 138  4.8397018 6.2831853}%
%
\special{pn 8}%
\special{ar 2064 1396 526 138  3.1415927 4.5848302}%
%
\special{pn 20}%
\special{pa 2616 1970}%
\special{pa 2612 1938}%
\special{pa 2610 1906}%
\special{pa 2606 1874}%
\special{pa 2604 1842}%
\special{pa 2600 1810}%
\special{pa 2598 1778}%
\special{pa 2596 1746}%
\special{pa 2596 1714}%
\special{pa 2594 1682}%
\special{pa 2594 1650}%
\special{pa 2594 1618}%
\special{pa 2594 1586}%
\special{pa 2594 1554}%
\special{pa 2594 1522}%
\special{pa 2594 1490}%
\special{pa 2596 1458}%
\special{pa 2596 1426}%
\special{pa 2598 1394}%
\special{pa 2600 1362}%
\special{pa 2604 1330}%
\special{pa 2606 1298}%
\special{pa 2610 1268}%
\special{pa 2612 1236}%
\special{pa 2616 1204}%
\special{pa 2622 1172}%
\special{pa 2624 1140}%
\special{pa 2630 1108}%
\special{pa 2636 1076}%
\special{pa 2640 1046}%
\special{pa 2646 1014}%
\special{pa 2652 982}%
\special{pa 2660 952}%
\special{pa 2666 920}%
\special{pa 2672 888}%
\special{pa 2680 858}%
\special{pa 2692 828}%
\special{pa 2698 796}%
\special{pa 2708 766}%
\special{pa 2716 736}%
\special{pa 2728 706}%
\special{pa 2740 674}%
\special{pa 2750 646}%
\special{pa 2758 622}%
\special{sp}%
%
\special{pn 8}%
\special{pa 2654 940}%
\special{pa 2638 968}%
\special{pa 2612 986}%
\special{pa 2582 998}%
\special{pa 2552 1006}%
\special{pa 2522 1016}%
\special{pa 2490 1022}%
\special{pa 2458 1026}%
\special{pa 2426 1030}%
\special{pa 2394 1034}%
\special{pa 2362 1036}%
\special{pa 2330 1038}%
\special{pa 2298 1038}%
\special{pa 2266 1040}%
\special{pa 2234 1040}%
\special{pa 2202 1038}%
\special{pa 2170 1038}%
\special{pa 2138 1038}%
\special{pa 2106 1036}%
\special{pa 2074 1034}%
\special{pa 2044 1032}%
\special{pa 2012 1028}%
\special{pa 1980 1024}%
\special{pa 1948 1020}%
\special{pa 1916 1014}%
\special{pa 1884 1010}%
\special{pa 1854 1004}%
\special{pa 1822 998}%
\special{pa 1790 990}%
\special{pa 1760 980}%
\special{pa 1730 972}%
\special{pa 1700 960}%
\special{pa 1670 946}%
\special{pa 1642 932}%
\special{pa 1618 910}%
\special{pa 1604 882}%
\special{pa 1604 878}%
\special{sp}%
%
\special{pn 8}%
\special{pa 2200 792}%
\special{pa 2232 798}%
\special{pa 2264 800}%
\special{pa 2296 804}%
\special{pa 2328 808}%
\special{pa 2360 814}%
\special{pa 2390 818}%
\special{pa 2422 824}%
\special{pa 2452 834}%
\special{pa 2484 842}%
\special{pa 2514 850}%
\special{pa 2544 862}%
\special{pa 2572 874}%
\special{pa 2600 890}%
\special{pa 2626 908}%
\special{pa 2646 934}%
\special{pa 2646 950}%
\special{sp}%
%
\special{pn 8}%
\special{pa 1594 882}%
\special{pa 1608 856}%
\special{pa 1632 836}%
\special{pa 1662 824}%
\special{pa 1690 812}%
\special{pa 1722 802}%
\special{pa 1752 794}%
\special{pa 1784 788}%
\special{pa 1816 786}%
\special{pa 1846 780}%
\special{pa 1878 778}%
\special{pa 1910 774}%
\special{pa 1942 774}%
\special{pa 1974 774}%
\special{pa 2006 772}%
\special{pa 2038 772}%
\special{pa 2064 772}%
\special{sp}%
%
\special{pn 8}%
\special{pa 1554 1870}%
\special{pa 1552 1838}%
\special{pa 1550 1806}%
\special{pa 1548 1774}%
\special{pa 1546 1742}%
\special{pa 1546 1710}%
\special{pa 1544 1678}%
\special{pa 1544 1646}%
\special{pa 1544 1614}%
\special{pa 1544 1582}%
\special{pa 1544 1550}%
\special{pa 1544 1518}%
\special{pa 1544 1486}%
\special{pa 1546 1454}%
\special{pa 1546 1422}%
\special{pa 1546 1390}%
\special{pa 1548 1358}%
\special{pa 1550 1326}%
\special{pa 1550 1294}%
\special{pa 1552 1262}%
\special{pa 1554 1230}%
\special{pa 1558 1198}%
\special{pa 1560 1166}%
\special{pa 1562 1134}%
\special{pa 1566 1102}%
\special{pa 1570 1070}%
\special{pa 1572 1040}%
\special{pa 1576 1008}%
\special{pa 1580 976}%
\special{pa 1584 944}%
\special{pa 1590 912}%
\special{pa 1594 880}%
\special{pa 1600 850}%
\special{pa 1606 818}%
\special{pa 1610 786}%
\special{pa 1618 754}%
\special{pa 1624 724}%
\special{pa 1630 692}%
\special{pa 1638 660}%
\special{pa 1646 630}%
\special{pa 1654 600}%
\special{pa 1662 568}%
\special{pa 1670 538}%
\special{pa 1674 528}%
\special{sp}%
%
\special{pn 8}%
\special{ar 2032 1448 492 1668  2.7534591 2.8488621}%
%
\special{pn 13}%
\special{ar 2994 1404 396 1732  2.6909842 2.7469348}%
%
\special{pn 8}%
\special{pa 2822 998}%
\special{pa 2428 1186}%
\special{dt 0.045}%
\special{sh 1}%
\special{pa 2428 1186}%
\special{pa 2496 1174}%
\special{pa 2476 1162}%
\special{pa 2480 1138}%
\special{pa 2428 1186}%
\special{fp}%
%
\special{pn 8}%
\special{pa 1078 1298}%
\special{pa 1352 1720}%
\special{dt 0.045}%
\special{sh 1}%
\special{pa 1352 1720}%
\special{pa 1332 1652}%
\special{pa 1324 1674}%
\special{pa 1300 1674}%
\special{pa 1352 1720}%
\special{fp}%
\put(28.7200,-10.2200){\makebox(0,0)[lb]{$t^F_r(B)$}}%
%
\special{pn 8}%
\special{ar 2702 1572 624 1480  3.6619480 3.8741101}%
\put(13.1000,-12.5700){\makebox(0,0)[rb]{$\exp^{\perp}(T^{\perp}_{p_0}B)$}}%
%
\special{pn 8}%
\special{pa 2904 1202}%
\special{pa 2622 1396}%
\special{dt 0.045}%
\special{sh 1}%
\special{pa 2622 1396}%
\special{pa 2688 1374}%
\special{pa 2666 1366}%
\special{pa 2666 1342}%
\special{pa 2622 1396}%
\special{fp}%
\put(29.5200,-12.4900){\makebox(0,0)[lb]{${\rm Exp}(\Phi\cdot{\widehat f}^F_{B,r}(v))$}}%
%
\special{pn 8}%
\special{pa 2040 456}%
\special{pa 2210 658}%
\special{dt 0.045}%
\special{sh 1}%
\special{pa 2210 658}%
\special{pa 2182 594}%
\special{pa 2176 618}%
\special{pa 2152 620}%
\special{pa 2210 658}%
\special{fp}%
\put(20.5500,-4.1600){\makebox(0,0)[rb]{$B$}}%
\put(20.2200,-15.8100){\makebox(0,0)[rt]{$p_0$}}%
%
\special{pn 8}%
\special{pa 3010 658}%
\special{pa 2702 764}%
\special{dt 0.045}%
\special{sh 1}%
\special{pa 2702 764}%
\special{pa 2772 760}%
\special{pa 2752 746}%
\special{pa 2758 724}%
\special{pa 2702 764}%
\special{fp}%
\put(30.5100,-7.3000){\makebox(0,0)[lb]{$\exp(P^{hol}_{f^F_{B,r}(v)})$}}%
%
\special{pn 8}%
\special{ar 3486 1412 1186 384  1.5502939 1.5655903}%
\special{ar 3486 1412 1186 384  1.6114794 1.6267757}%
\special{ar 3486 1412 1186 384  1.6726648 1.6879612}%
\special{ar 3486 1412 1186 384  1.7338503 1.7491467}%
\special{ar 3486 1412 1186 384  1.7950358 1.8103321}%
\special{ar 3486 1412 1186 384  1.8562212 1.8715176}%
\special{ar 3486 1412 1186 384  1.9174067 1.9327031}%
\special{ar 3486 1412 1186 384  1.9785922 1.9938885}%
\special{ar 3486 1412 1186 384  2.0397776 2.0550740}%
\special{ar 3486 1412 1186 384  2.1009631 2.1162595}%
\special{ar 3486 1412 1186 384  2.1621486 2.1774449}%
\special{ar 3486 1412 1186 384  2.2233340 2.2386304}%
\special{ar 3486 1412 1186 384  2.2845195 2.2998159}%
\special{ar 3486 1412 1186 384  2.3457050 2.3610013}%
\special{ar 3486 1412 1186 384  2.4068904 2.4221868}%
\special{ar 3486 1412 1186 384  2.4680759 2.4833723}%
\special{ar 3486 1412 1186 384  2.5292614 2.5445578}%
\special{ar 3486 1412 1186 384  2.5904469 2.6057432}%
\special{ar 3486 1412 1186 384  2.6516323 2.6669287}%
\special{ar 3486 1412 1186 384  2.7128178 2.7281142}%
\special{ar 3486 1412 1186 384  2.7740033 2.7892996}%
\special{ar 3486 1412 1186 384  2.8351887 2.8504851}%
\special{ar 3486 1412 1186 384  2.8963742 2.8980492}%
%
\special{pn 8}%
\special{pa 3524 1796}%
\special{pa 3574 1796}%
\special{dt 0.045}%
\special{sh 1}%
\special{pa 3574 1796}%
\special{pa 3508 1776}%
\special{pa 3522 1796}%
\special{pa 3508 1816}%
\special{pa 3574 1796}%
\special{fp}%
\put(29.0400,-18.4000){\makebox(0,0)[lt]{in fact}}%
%
\special{pn 8}%
\special{ar 4276 2004 608 576  0.0000000 6.2831853}%
%
\special{pn 8}%
\special{ar 4276 2014 600 142  0.0000000 6.2831853}%
%
\special{pn 8}%
\special{ar 1620 2520 3912 1206  4.8409042 4.9492437}%
%
\special{pn 8}%
\special{ar 1620 2520 3930 1206  4.7280203 4.8241104}%
%
\special{pn 8}%
\special{ar 1594 2520 3816 1206  4.6801851 4.7108715}%
%
\special{pn 8}%
\special{ar 2984 1590 906 1668  2.8105390 2.9150004}%
%
\special{pn 8}%
\special{sh 0.300}%
\special{ar 2596 1404 28 24  0.0000000 6.2831853}%
%
\special{pn 20}%
\special{ar 4248 2022 432 712  0.1807147 0.7492965}%
%
\special{pn 20}%
\special{ar 4398 2058 274 690  5.3880514 6.2831853}%
\special{ar 4398 2058 274 690  0.0000000 0.0558159}%
%
\special{pn 20}%
\special{ar 4682 2162 252 418  2.1496249 3.1415927}%
%
\special{pn 20}%
\special{ar 4752 2156 324 810  3.1499085 4.0823858}%
%
\special{pn 8}%
\special{sh 0.300}%
\special{ar 2606 1404 36 32  0.0000000 6.2831853}%
\put(41.1000,-10.2000){\makebox(0,0)[lb]{${\rm Exp}(\Phi\cdot{\widehat f}^F_{B,r}(v))$}}%
\put(39.9200,-27.3700){\makebox(0,0)[lt]{$f^F_{B,r}(\pi_B^{-1}(p_0))$}}%
%
\special{pn 8}%
\special{pa 4310 1060}%
\special{pa 4530 1540}%
\special{dt 0.045}%
\special{sh 1}%
\special{pa 4530 1540}%
\special{pa 4520 1472}%
\special{pa 4508 1492}%
\special{pa 4484 1488}%
\special{pa 4530 1540}%
\special{fp}%
\end{picture}%
\hspace{7truecm}}

\vspace{0.75truecm}

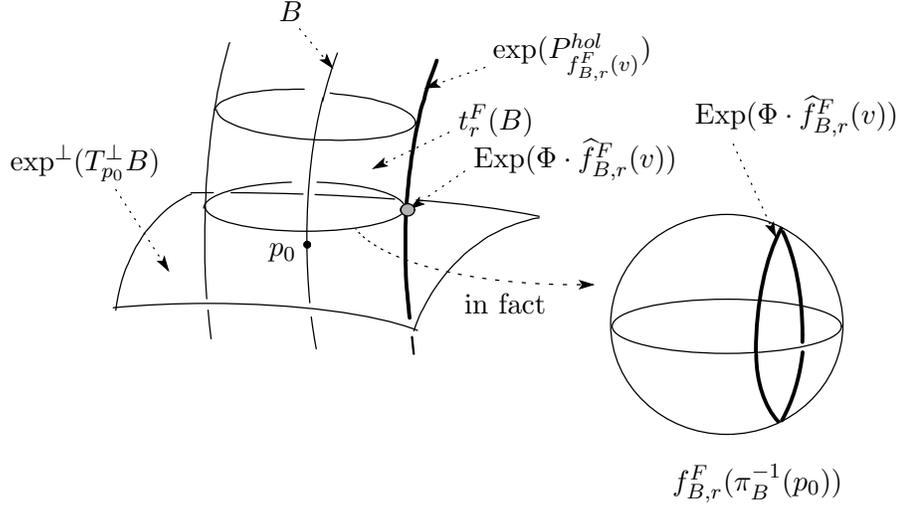
\centerline{{\bf Figure 5.}}

\vspace{0.5truecm}

At the end of this section, we give the list of all (commutative) Hermann actions $H\curvearrowright G/K$ 
of chomogeneity one on irreducible symmetric spaces $G/K$ of non-compact type and rank greater than one 
and the reflective singular orbit $B=H(eK)$ (see [BT,Koi3]).  

\vspace{0.5truecm}

$$\begin{tabular}{|c|c|c|}
\hline
{\scriptsize$H\curvearrowright G/K$} & {\scriptsize$B=H(eK)$} & {\scriptsize Remark}\\
\hline
{\scriptsize$SO_0(k-1,n-k)\curvearrowright$} 
& {\scriptsize$SO_0(k-1,n-k)/SO(k-1)\times SO(n-k)$} & {\scriptsize $2\leq k<n/2$}\\
{\scriptsize$SO_0(k,n-k)/SO(k)\times SO(n-k)$} & & 
\\
\hline
{\scriptsize$SO_0(k,n-k-1)\curvearrowright$} 
& {\scriptsize$SO_0(k,n-k-1)/\times SO(k)\times SO(n-k-1)$} & 
{\scriptsize $2\leq k<n/2$}\\
{\scriptsize$SO_0(k,n-k)/SO(k)\times SO(n-k)$} & & 
\\
\hline
{\scriptsize$SO_0(k-1,k)\curvearrowright$} 
& {\scriptsize$SO_0(k-1,k)/SO(k-1)\times SO(k)$} & 
{\scriptsize $k\geq 3$}\\
{\scriptsize$SO_0(k,k)/SO(k)\times SO(k)$} && \\
\hline
{\scriptsize$SL(3,{\Bbb R})\cdot U(1)\curvearrowright$} 
& {\scriptsize$(SL(3,{\Bbb R})/SO(3))\times{\Bbb R}$} & \\
{\scriptsize$SO_0(3,3)/SO(3)\times SO(3)$} &&\\
\hline
{\scriptsize$SU(1,n-1)\cdot U(1)\curvearrowright$} 
& {\scriptsize$SU(1,n-1)/S(U(1)\times U(n-1))$} & 
{\scriptsize $n\geq 3$}\\
{\scriptsize$SO_0(2,2n-2)/SO(2)\times SO(2n-2)$} &&\\
\hline
{\scriptsize$SU(k-1,n-k)\curvearrowright$} 
& {\scriptsize$SU(k-1,n-k)/S(U(k-1)\times U(n-k))$} & {\scriptsize $2\leq k<n/2$}\\
{\scriptsize$SU(k,n-k)/S(U(k)\times U(n-k))$} && 
\\
\hline
{\scriptsize$SU(k,n-k-1)\curvearrowright$} 
& {\scriptsize$SU(k,n-k-1)/S(U(k)\times U(n-k-1))$} & {\scriptsize $2\leq k<n/2$}\\
{\scriptsize$SU(k,n-k)/S(U(k)\times U(n-k))$} && 
\\
\hline
{\scriptsize$SU(k-1,k)\curvearrowright$} 
& {\scriptsize$SU(k-1,k)/S(U(k-1)\times U(k))$} & 
{\scriptsize $k\geq 3$}\\
{\scriptsize$SU(k,k)/S(U(k)\times U(k))$} &&\\
\hline
{\scriptsize$Sp(1,n-1)\curvearrowright$} 
& {\scriptsize$Sp(1,n-1)/Sp(1)\times Sp(n-1)$} & {\scriptsize $n\geq 3$}\\
{\scriptsize$SU(2,2n-2)/S(U(2)\times U(2n-2))$} &&\\
\hline
\end{tabular}$$

\vspace{0.3truecm}

\centerline{{\bf Table 1.}}

\vspace{0.5truecm}

$$\begin{tabular}{|c|c|c|}
\hline
{\scriptsize$H\curvearrowright G/K$} & {\scriptsize$B=H(eK)$} & {\scriptsize Remark}\\
\hline
{\scriptsize$Sp(k-1,n-k)\curvearrowright$} 
& {\scriptsize$Sp(k-1,n-k)/Sp(k-1)\times Sp(n-k)$} & {\scriptsize $2\leq k<n/2$}\\
{\scriptsize$Sp(k,n-k)/Sp(k)\times Sp(n-k)$} &&\\
\hline
{\scriptsize$Sp(k,n-k-1)\curvearrowright$} 
& {\scriptsize$Sp(k,n-k-1)/Sp(k)\times Sp(n-k-1)$} & {\scriptsize $2\leq k<n/2$}\\
{\scriptsize$Sp(k,n-k)/Sp(k)\times Sp(n-k)$} &&\\
\hline
{\scriptsize$Sp(k-1,k)\curvearrowright$} 
& {\scriptsize$Sp(k-1,k)/Sp(k-1)\times Sp(k)$} & 
{\scriptsize $k\geq 2$}\\
{\scriptsize$Sp(k,k)/Sp(k)\times Sp(k)$} &&\\
\hline
{\scriptsize$Sp(2,{\Bbb C})\curvearrowright$} 
& {\scriptsize$Sp(2,{\Bbb C})/Sp(2)$} & \\
{\scriptsize$Sp(2,2)/Sp(2)\times Sp(2)$} &&\\
\hline
{\scriptsize$SL(n-1,{\Bbb R})\cdot{\Bbb R}_{\ast}\curvearrowright SL(n,{\Bbb R})/SO(n)$} & 
{\scriptsize$(SL(n-1,{\Bbb R})/SO(n-1))\times{\Bbb R}$} & {\scriptsize$n\geq 3$}\\
\hline
{\scriptsize$SU^{\ast}(2n-2)\cdot{\Bbb R}_{\ast}\curvearrowright SU^{\ast}(2n)/Sp(n)$} 
& {\scriptsize$(SU^{\ast}(2n-2)/Sp(n-1))\times{\Bbb R}$} & {\scriptsize $n\geq 3$}\\
\hline
{\scriptsize$SO^{\ast}(2n-2)\curvearrowright SO^{\ast}(2n)/U(n)$} 
& {\scriptsize$SO^{\ast}(2n-2)/U(n-1)$} & {\scriptsize $n\geq 4$}\\
\hline
{\scriptsize$SU(1,3)\cdot U(1)\curvearrowright SO^{\ast}(8)/U(4)$} 
& {\scriptsize$SU(1,3)/S(U(1)\times U(3))$} & \\
\hline
{\scriptsize$SU^{\ast}(4)\cdot U(1)\curvearrowright SO^{\ast}(8)/U(4)$} 
& {\scriptsize$SU^{\ast}(4)/Sp(2)$} & \\
\hline
{\scriptsize$Sp(1,{\Bbb R})\times Sp(n-1,{\Bbb R})\curvearrowright$} 
& {\scriptsize$(Sp(1,{\Bbb R})/U(1))\times(Sp(n-1,{\Bbb R})/U(n-1))$} & {\scriptsize $n\geq 3$}\\
{\scriptsize$Sp(n,{\Bbb R})/U(n)$}&&\\
\hline
{\scriptsize$SO(n-1,{\Bbb C})\curvearrowright SO(n,{\Bbb C})/SO(n)$} 
& {\scriptsize$SO(n-1,{\Bbb C})/SO(n-1)$} & {\scriptsize $n\geq 5$}\\
\hline
{\scriptsize$SL(3,{\Bbb C})\cdot SO(2,{\Bbb C})\curvearrowright SO(6,{\Bbb C})/SO(6)$} 
& {\scriptsize$(SL(3,{\Bbb C})/SU(3))\times(SO(2,{\Bbb C})/SO(2))$} & \\
\hline
{\scriptsize$SL(n-1,{\Bbb C})\times{\Bbb C}_{\ast}\curvearrowright$} 
& {\scriptsize$(SL(n-1,{\Bbb C})/SU(n-1))\times{\Bbb R}$} & {\scriptsize $n\geq 3$}\\
{\scriptsize$SL(n,{\Bbb C})/SU(n)$} &&\\
\hline
{\scriptsize$SL(3,{\Bbb R})\curvearrowright SL(3,{\Bbb C})/SU(3)$} 
& {\scriptsize$SL(3,{\Bbb R})/SO(3)$} & \\
\hline
{\scriptsize$Sp(2,{\Bbb C})\curvearrowright SL(4,{\Bbb C})/SU(4)$} 
& {\scriptsize$Sp(2,{\Bbb C})/Sp(2)$} & \\
\hline
{\scriptsize$SU(1,2)\curvearrowright SL(3,{\Bbb C})/SU(3)$} 
& {\scriptsize$SU(1,2)/S(U(1)\times U(2))$} & \\
\hline
{\scriptsize$Sp(1,{\Bbb C})\times Sp(n-1,{\Bbb C})\curvearrowright$} 
& {\scriptsize$(Sp(1,{\Bbb C})/Sp(1))\times(Sp(n-1,{\Bbb C})/Sp(n-1))$} & {\scriptsize $n\geq 3$}\\
{\scriptsize$Sp(n,{\Bbb C})/Sp(n)$} &&\\
\hline
{\scriptsize$F_4^4\curvearrowright E_6^2/SU(6)\cdot SU(2)$} 
& {\scriptsize$F_4^4/Sp(3)\cdot Sp(1)$} & \\
\hline
{\scriptsize$F_4^{-20}\curvearrowright E_6^{-14}/Spin(10)\cdot U(1)$} 
& {\scriptsize$F_4^{-20}/Spin(9)$} & \\
\hline
{\scriptsize$SU^{\ast}(6)\cdot SU(2)\curvearrowright E_6^{-26}/F_4$} 
& {\scriptsize$SU^{\ast}(6)\cdot SU(2)/Sp(3)\cdot Sp(1)$} & \\
\hline
{\scriptsize$SO_0(1,9)\cdot U(1)\curvearrowright E_6^{-26}/F_4$} 
& {\scriptsize$SO_0(1,9)\cdot U(1)/SO(1)\times SO(9)$} & \\
\hline
{\scriptsize$SO_0(4,5)\curvearrowright F_4^4/Sp(3)\cdot Sp(1)$} 
& {\scriptsize$SO_0(4,5)/SO(4)\times SO(5)$} & \\
\hline
{\scriptsize$SO(9,{\Bbb C})\curvearrowright F_4^{\Bbb C}/F_4$} 
& {\scriptsize$SO(9,{\Bbb C})/SO(9)$} & \\
\hline
\end{tabular}$$

\vspace{0.5truecm}

\centerline{{\bf Table 1(continued).}}

\vspace{0.5truecm}

The dual action $H^{\ast}\curvearrowright G^{\ast}/K$ of a (commutative) Hermann action 
$H\curvearrowright G/K$ is defined in a natural manner, where $G^{\ast}$ is the compact dual of $G$ 
with respect to $K$ and $H^{\ast}$ is the compact dual of $H$ with respect to $H\cap K$.  
The dual actions of Hermann actions in Table 1 are all of (commutative) Hermann actions 
of cohomogeneity one on irreducible symmetric space of compact type and rank greater than one.  

\section{The equivalenceness of the anisotropic equifocality and the anisotropic 
isoparametricness} 
We use the notations in Sections 1-5.  
Assume that $\partial\, M=\emptyset$.  
In the case where $\widetilde M$ is a Euclidean space, J. Ge and H. Ma ([GM]) showed that 
$f$ is of constant anisotropic principal curvatures (i.e., 
anisotropic equifocal) if and only if $f$ is anisotropic isoparametric.  
We obtain the following similar result in the case where $\widetilde M$ is 
a symmetric space of non-negative curvature.   


\vspace{0.5truecm}

\noindent
{\bf Theorem 6.1.} {\sl In the case where $\widetilde M$ is a symmetric space of non-negative 
curvature, the anisotropic equifocality is equivalent to the anisotropic isoparametricness.}

\vspace{0.5truecm}

Assume that $\widetilde M$ is a symmetric space $G/K$ of non-negative curvature, 
where $G$ and $K$ are as in the previous section.  
We shall prove this theorem by reducing to the investigation of the lift of $f(M)$ by 
a Riemannian submersion of a Hilbert space onto $G/K$.  
We suffice to show the statement of this theorem in the case where $G/K$ is simply connected.  
In the sequel, we assume that $G/K$ is simply connected.  
In this case, $G/K$ is decomposed irreducibly as 
$\displaystyle{G/K=(\mathop{\Pi}_{i=1}^{\it l}G_i/K_i)\times{\Bbb R}^r}$ 
($G_i/K_i\,:\,$ a simply connected irreducible symmetric space of compact type, $r\,:\,$ 
a non-negative integer).  
Let $H^0([0,1],\mathfrak g_i)$ be the (separable) Hilbert space of all 
$L^2$-integrable paths in the Lie algebra $\mathfrak g_i$ of $G_i$ 
(having $[0,1]$ as the domain) and 
$H^1([0,1],G_i)$ the Hilbert Lie group of all $H^1$-paths in $G_i$ 
(having $[0,1]$ as the domain), where we give $G_i$ the bi-invariant metric inducing the metric 
of $G_i/K_i$ and $\mathfrak g_i$ the ${\rm Ad}(G)$-invariant inner product compatible with the 
bi-invariant metric.  
Let $\phi_i:H^0([0,1],\mathfrak g_i)\to G_i$ be the parallel transport map for 
$G_i$, that is, $\phi(u):=g_u(1)\,\,(u\in H^0([0,1],\mathfrak g_i))$, 
where $g_u$ is the element of $H^1([0,1],G_i)$ with $g_u(0)=e_i$ and 
$(R_{g_u(t)})_{\ast}^{-1}(g_u'(t))=u(t)\,\,(0\leq t\leq 1)$, where $e_i$ is 
the identity element of $G_i$.  
See [TT], [Ch], [HLO] and [Koi1,3] about the investigation of submanifold geometry in a symmetric 
space of compact type by using the parallel transport map.  
The group $H^1([0,1],G_i)$ acts on $H^0([0,1],\mathfrak g_i)$ 
isometrically as the action of the Gauge transformations to connections as follows:
$$\begin{array}{c}
\displaystyle{(g\ast u)(t):={\rm Ad}(g(t))(u(t))-(R_{g(t)})_{\ast}^{-1}(g'(t))}\\
\displaystyle{(g\in H^1([0,1],G_i),\,u\in H^0([0,1],\mathfrak g_i),\,\,0\leq t\leq1).}
\end{array}$$
Also, let $\pi_i:G_i\to G_i/K_i$ be the natural projection.  
Set $\widehat{\phi}_i:=\pi_i\circ\phi_i$, which is a Riemannian submersion of 
$H^0([0,1],\mathfrak g_i)$ onto $G_i/K_i$.  
Set $P(G_i,e_i\times K_i):=\{g\in H^1([0,1],G_i)\,\vert(g(0),g(1))\in \{e_i\}\times K_i\}$.  
This subgroup $P(G_i,e_i\times K_i)$ acts on $H^0([0,1],\mathfrak g_i)$ freely and 
$\widehat{\phi}_i:H^0([0,1],\mathfrak g_i)\to G_i/K_i$ is regarded as 
a $P(G_i,e_i\times K_i)$-bundle.  
For simplicity, set $\displaystyle{V:=(\mathop{\Pi}_{i=1}^{\it l}H^0([0,1],
\mathfrak g_i))\times{\Bbb R}^r}$, 
$\displaystyle{\widehat G:=\mathop{\Pi}_{i=1}^{\it l}P(G_i,e_i\times K_i)}$ 
and $\displaystyle{\widehat{\phi}:=(\mathop{\Pi}_{i=1}^{\it l}
\widehat{\phi}_i)\times{\rm id}_{{\Bbb R}^r}}$.  
Note that $\widehat{\phi}:V\to G/K$ is regarded as a $\widehat G$-bundle.  
We consider $\widehat F:=\widetilde F\circ\widehat{\phi}_{\ast}:TV\to{\Bbb R}$ 
as a parametric Lagrangian of $V$.  
Note that $\widehat F$ is not holonomy invariant.  
Denote by $\widehat M$ the induced bundle 
$f^{\ast}V:=\{(x,u)\in M\times V\,\vert\,f(x)=\widehat{\phi}(u)\}$ of this bundle 
$\widehat{\phi}:V\to G/K$ by $f$ and define an immersion 
$\widehat f:\widehat M\hookrightarrow H^0([0,1],\mathfrak g)$ by 
$\widehat f(x,u):=u\,\,((x,u)\in \widehat M)$.  
Note that $\widehat f(\widehat M)={\widehat{\phi}}^{-1}(f(M))$.  
In 1989, Terng ([T]) introduced the notion of a proper Fredholm submanifold in the Hilbert space 
and, in 2006, Heintze-Liu-Olmos ([HLO]) introduced the notion of a regularizable submanifold 
in the Hilbert space.  
According to Lemma 6.2 of [HLO], the hypersurface $\widehat f:\widehat M\hookrightarrow V$ is 
a regularizable (proper Fredholm) hypersurface.  The horizontal lift $\xi^L$ of $\xi$ is a 
unit normal vector field of $\widehat f$.  
Denote by $\widehat A$ the shape operator of $\widehat f$ for $\xi^L$.  
Since $\widehat f:\widehat M\hookrightarrow V$ is 
a regularizable hypersurface, for any $(x,u)\in\widehat M$, 
$\widehat A_{(x,u)}$ is a compact self-adjoint operator, and 
the regularized trace ${\rm Tr}_r\,\widehat A_{(x,u)}$ of 
$\widehat A_{(x,u)}$ and the trace ${\rm Tr}\,{\widehat A}_{(x,u)}^2$ 
exist.  See [HLO] (or [Koi4]) about the definition of the regularized trace.  
The regularized mean curvature $\widehat H$ of $\widehat f$ is defined by 
$\widehat H_{(x,u)}:={\rm Tr}_r\,\widehat A_{(x,u)}\,\,((x,u)\in\widehat M)$.  
Denote by $\widehat{\phi}_{\widehat M}$ the natural projection of $\widehat M$ onto $M$ 
(i.e., $\widehat{\phi}_{\widehat M}(x,u)=x\,\,((x,u)\in\widehat M))$.  
Define a transversal vector field $\widehat{\xi}_F$ of $\widehat f$ by 
$$(\widehat{\xi}_F)_{(x,u)}:=(F\circ\nu\circ\widehat{\phi}_{\widehat M})(x,u)\xi^L_u
+\widehat f_{\ast}({\rm grad}(F\circ\nu\circ\widehat{\phi}_{\widehat M}))\quad\,
((x,u)\in\widehat M).\leqno{(6.1)}$$
We call $\widehat{\xi}_F$ a {\it anisotropic transversal vector field} of $\widehat f$.  
Note that $\widehat F(\xi_u^L)=(F\circ\nu\circ\widehat{\phi}_{\widehat M})(x,u)$.  
Denote by $\widehat{\nabla}$ the Riemannian connection of $V$ and 
$\bar{\nabla}$ the Riemannian connection of the metric of $\widehat M$ induced by $\widehat f$.  
For any $X\in T_{(x,u)}\widehat M$, we can show 
$${\widehat{\nabla}}^{\widehat f}_X\widehat{\xi}_F
\equiv -(F\circ\nu\circ\widehat{\phi}_{\widehat M})(x,u)\widehat f_{\ast}
(\widehat A_{(x,u)}X)+\widehat f_{\ast}\left(\bar{\nabla}_X{\rm grad}
(F\circ\nu\circ\widehat{\phi}_{\widehat M})\right)\quad\,\,({\rm mod}\,T^{\perp}\widehat M),$$
where ${\widehat{\nabla}}^{\widehat f}$ is the covariant derivative along $\widehat f$ for 
$\widehat{\nabla}$.  
So, define a $(1,1)$-tensor field ${\widehat A}^F$ on 
$\widehat M$ by 
$${\widehat A}^F_{(x,u)}X
=(F\circ\nu\circ\widehat{\phi}_{\widehat M})(x,u)\widehat A_{(x,u)}X
-\bar{\nabla}_X{\rm grad}(F\circ\nu\circ\widehat{\phi}_{\widehat M})\leqno{(6.2)}$$
for any $(x,u)\in \widehat M$ and $X\in T_{(x,u)}\widehat M$.  
We call ${\widehat A}^F$ an {\it anisotropic shape operator} of 
$\widehat f$ and the eigenvalues of ${\widehat A}^F$ {\it anisotropic 
principal curvatures} of $\widehat f$.  
Since $F\circ\nu\circ\widehat{\phi}_{\widehat M}$ is 
$\widehat G$-invariant, $\bar{\nabla}{\rm grad}
(F\circ\nu\circ\widehat{\phi}_{\widehat M})$ is a compact self-adjoint regularizable operator.  
On the other hand, so is also 
${\widehat A}_{(x,u)}$.  Hence so is also ${\widehat A}^F$.  
Define a function $\widehat H_F$ over $\widehat M$ by 
$\widehat H_F:={\rm Tr}_r\,{\widehat A}^F$, 
where ${\rm Tr}_r\,{\widehat A}^F$ is the regularized trace of 
${\widehat A}^F$.  We call $\widehat H_F$ an 
{\it anisotropic regularized mean curvature of} $\widehat f$.  
Define a map $\widehat f_t:\widehat M\to V$ by 
$$\widehat f_t(x,u):=\widehat f(x,u)+t(\widehat{\xi}_F)_{(x,u)}
\quad\,\,((x,u)\in\widehat M).$$
Since $\widehat f(\widehat M)$ is $\widehat G$-invariant and 
$\widehat{\phi}(\widehat f(\widehat M))$ is compact, 
it is shown that $\widehat f_t$ is an immersion for each $t$ sufficiently close to zero.  
If $\widehat f_t$ is an immersion, then we call 
$\widehat f_t:\widehat M\hookrightarrow V$ an 
{\it anisotropic parallel hypersurface} of 
$\widehat f:\widehat M\hookrightarrow V$.  
If, for each $t$ sufficiently close to zero, 
$\widehat f_t:\widehat M\hookrightarrow V$ is 
of constant anisotropic regularized mean curvature, we call 
$\widehat f:\widehat M\hookrightarrow V$ an 
{\it anisotropic isoparametric hypersurface}.  
Also, if the spectrum of ${\widehat A}^F_{(x,u)}$ is 
independent of the choice of $(x,u)\in\widehat M$, then 
we call $\widehat f:\widehat M\hookrightarrow V$ a 
{\it hypersurface with constant} {\it anisotropic principal curvatures}.  
The following facts follow directly.  

\vspace{0.5truecm}

\noindent
{\bf Lemma 6.2.} {\sl 
{\rm(i)} $\widehat{\xi}_F$ is the horizontal lift of $\xi_F$.  

{\rm(ii)} $\widehat H_F=H_F\circ\widehat{\phi}_{\widehat M}$ holds.}

\vspace{0.5truecm}

\noindent
{\it Proof.} The statement (i) follows from $(2.1)$ and $(6.1)$ directly.  
Also, the statement (ii) follows from $(2.3)$ and $(6.2)$ directly.  
\hspace{7.7truecm}q.e.d.

\vspace{0.5truecm}

By using this lemma, we shall prove Thoerem 6.1.  

\vspace{0.5truecm}

\noindent
{\it Proof of Theorem 6.1.} First we shall show that $f$ is anisotropic 
equifocal if and only if $\widehat f$ is of constant 
anisotropic principal curvatures.  
Denote by $\gamma^F_{(x,u)}$ the geodesic in 
$V$ whose initial velocity vector is equal to 
$(\widehat{\xi}_F)_{(x,u)}$.  
Take $X\in T_{(x,u)}\widehat M$.  The anisotropic $\widehat M$-Jacobi field $Y_X$ along 
$\gamma^F_{(x,u)}$ with $Y_X(0)={\widehat f}_{\ast}X$ and 
$Y_X'(0)=-{\widehat f}_{\ast}({\widehat A}^F_{(x,u)}X)$ is described as 
$$Y_X(s)={\widehat f}_{\ast}X-s{\widehat f}_{\ast}({\widehat A}^F_{(x,u)}(X)).\leqno{(6.3)}$$
From this description, anisotropic focal radii of 
$\widehat f$ at $(x,u)$ are equal to the inverse numbers of anisotropic principal curvatures.  
Also, since $\widehat{\phi}$ is a Riemannian submersion and $\widehat{\xi}_F$ is the horizontal 
lift of $\xi_F$, the anisotropic focal radii of $f$ at $x$ are equal to those of $\widehat f$ 
with considering their multiplicities.  
From these facts, it follows that $f$ is anisotropic equifocal if and only if $\widehat f$ has 
constant anisotropic principal curvatures.  

Next we shall show that $\widehat f$ has constant anisotropic 
principal curvatures if and only if $\widehat f$ is 
anisotropic isoparametric.  
Take a positive number $\varepsilon$ such that $\widehat f_t\,\,
(-\varepsilon<t<\varepsilon)$ are immersions.  
Denote by $\widehat{\xi}_F^t$ and ${\widehat A}^F_t$ 
the anisotropic transversal vector field and 
the anisotropic shape operator of the parallel hypersurface 
$\widehat f_t:\widehat M\hookrightarrow V$ 
($-\varepsilon<t<\varepsilon$), respectively.  
Also, denote by ${\widehat H}_F^t$ the 
anisotropic regularized mean curvature of $\widehat f_t$.  
Easily we can show that $\widehat{\xi}_F^t=\widehat{\xi}_F$ in $V$.  
Assume that ${\widehat A}^F_{(x,u)}(X)=\lambda X$.  
Let $Y_X$ be the anisotropic $\widehat M$-Jacobi field along 
$\gamma^F_{(x,u)}$ with 
$Y_X(0)={\widehat f}_{\ast}X$ and $Y_X'(0)=-{\widehat f}_{\ast}({\widehat A}^F_{(x,u)}X)$.  
Since $Y_X(t)=({\widehat f}_t)_{\ast}X,\,Y'_X(t)=-({\widehat f}_t)_{\ast}
((\widehat A^F_t)_{(x,u)}(X))$, and since $Y_X$ is described as in $(6.3)$, 
we have 
$$({\widehat f}_t)_{\ast}X=(1-t\lambda){\widehat f}_{\ast}X$$
and 
$$({\widehat f}_t)_{\ast}((\widehat A^F_t)_{(x,u)}(X))=\lambda{\widehat f}_{\ast}X.$$
From these relation, we obtain 
$$({\widehat A}^F_t)_{(x,u)}(X)
=\frac{\lambda}{1-t\lambda}X.\leqno{(6.4)}$$
Denote by ${\rm Spec}(\cdot)$ the spectrum of $(\cdot)$.  
Set $m_{\lambda}:={\rm dim}\,{\rm Ker}({\widehat A}^F_{(x,u)}-\lambda\,{\rm id})$ for each 
$\lambda\in{\rm Spec}\,{\widehat A}^F_{(x,u)}$.  
From $(6.4)$, we have 
$${\rm Spec}\,({\widehat A}^F_t)_{(x,u)}
=\left\{\left.\frac{\lambda}{1-t\lambda}\,\right\vert\,\lambda\in{\rm Spec}\,
{\widehat A}^F_{(x,u)}\right\}$$
and 
$\displaystyle{{\rm dim}\,{\rm Ker}\left(({\widehat A}^F_t)_{(x,u)}-\frac{\lambda}{1-t\lambda}\,{\rm id}
\right)=m_{\lambda}}$.  
Hence we have 
$$({\widehat H}_F^t)_{(x,u)}
=\sum_{\lambda\in{\rm Spec}\,{\widehat A}^F_{(x,u)}}
\frac{m_{\lambda}\lambda}{1-t\lambda},$$
where the right-hand side means the regularized series.  
Hence, by using Lemma 4.1 in [HLO], we can show that 
$\widehat f$ has constant anisotropic principal curvatures if and only 
if $\widehat f$ is anisotropic isoparametric.  

Finally we shall show that $f$ is anisotropic 
isoparametric if and only if so is $\widehat f$.  
Denote by $H_F^t$ the anisotropic mean curvature of $f_t$.  
According to (i) of Lemma 6.2, we have 
$\widehat{\phi}\circ\widehat f_t
=f_t\circ{\widehat{\phi}}_{\widehat M}$.  
Hence we have $\widehat H_F^t=H_F^t\circ\widehat{\phi}_{\widehat M}$ in similar to (ii) of 
Lemma 6.2.  This implies that $f$ is anisotropic isoparametric if and only if 
so is $\widehat f$.  Therefore the statement of Theorem 6.1 follows.  
\hspace{10.4truecm}q.e.d.




\vspace{1truecm}

\centerline{{\bf References}}

\vspace{1truecm}

{\small

%





\noindent
[BCO] J. Berndt, S. Console and C. Olmos,
Submanifolds and holonomy, Research Notes 
in Mathe-

matics 434, CHAPMAN $\&$ HALL/CRC Press, Boca Raton, London, New York 
Washington, 

2003.



\noindent
[BT] J. Berndt and H. Tamaru, Cohomogeneity one actions on noncompact 
symmetric spaces 

with a totally geodesic singular orbit, 
Tohoku Math. J. {\bf 56} (2004) 163-177.

\noindent
[Ch] U. Christ, 
Homogeneity of equifocal submanifolds, J. Differential Geom. {\bf 62} (2002) 1-15.

\noindent
[Cl] U. Clarenz,
Enclosure theorems for extremals of elliptic parametric 
functionals, Calc. Var. {\bf 15} 

(2002) 313-324. 

\noindent
[CW] S. Carter and A. West, 
Partial tubes about immersed manifolds, Geom. Dedicata {\bf 54} (1995) 

145-169. 

%
%
\noindent
[F1] H. Federer, 
Geometric Measure Theory, Grundlehren math. Wiss. {\bf 153}. Berlin, Heidelberg, 

New York:Springer, 1969.

\noindent
[F2] H. Federer, 
Real flat chains, cochains and variational problems, Indiana Univ. Math. J. {\bf 24} 

(1974) 351-407.

\noindent
[F3] H. Federer, 
Colloquium Lectures on Geometric Measure Theory, Bulletin of the Amer. Math. 

Soc. {\bf 84(3)}, May 1978.

\noindent
[GM] J. Ge and H. Ma, 
Anisotropic isoparametric hypersurfaces in Euclidean spaces, Ann. Glob. 

Anal. Geom. {\bf 41} (2012) 347-355.



%

\noindent
[HLO] E. Heintze, X. Liu and C. Olmos, Isoparametric submanifolds and a 
Chevalley type rest-

riction theorem, Integrable systems, geometry, and topology, 151-190, 
AMS/IP Stud. Adv. 

Math. 36, Amer. Math. Soc., Providence, RI, 2006.

\noindent
[H] S. Helgason, 
Differential geometry, Lie groups and symmetric spaces, 
Academic Press, New 

York, 1978.



\noindent
[Koi1] N. Koike,
On proper Fredholm submanifolds in a Hilbert space 
arising from submanifolds 

in a symmetric space, 
Japan. J. Math. {\bf 28} (2002) 61--80.

\noindent
[Koi2] N. Koike, 
Actions of Hermann type and proper complex equifocal 
submanifolds, Osaka J. 

Math. {\bf 42} (2005) 599-611.

\noindent
[Koi3] N. Koike, 
Complex hyperpolar actions with a totally geodesic orbit, 
Osaka J. Math. {\bf 44} 

(2007) 491-503.

\noindent
[Koi4] N. Koike, 
Collapse of the mean curvature flow for equifocal 
submanifolds, Asian J. Math. 

{\bf 15} (2011) 101-128.

%

\noindent
[Kol] A. Kollross, A classification of hyperpolar and cohomogeneity one 
actions, Trans. Amer. 

Math. Soc. {\bf 354} (2002) 571-612.

\noindent
[KP1] M. Koiso and B. Palmer, 
Geometry and stability of surfaces with constant anisotropic 

mean curvature, Indiana Univ. Math. J. {\bf 54} (2005) 1817-1852.

\noindent
[KP2] M. Koiso and B. Palmer, 
Anisotropic capillary surfaces with wetting energy, Calc. Var. {\bf 29} 

(2007) 295-345.

\noindent
[KP3] M. Koiso and B. Palmer, 
Anisotropic umbilic points and Hopf's Theorem for surfaces with 

constant anisotropic mean curvature, Indiana Univ. Math. J. {\bf 59} (2010) 
79-90.

\newpage

\noindent
[LM] J.H.S. de Lira and M. Melo, 
Hypersurfaces with constant anisotropic mean curvature in Rie-

mannian manifolds, Calc. Var. {\bf 50} (2014) 335-364.




\noindent
[Pala] R. S. Palais,
Morse theory on Hilbert manifolds, 
Topology {\bf 2} (1963) 299--340.

\noindent
[Palm] B. Palmer, 
Stability of the Wulff shape, Proc. Amer. Math. Soc. {\bf 126} 
(1998) 3661-3667.

\noindent
[PT] R.S. Palais and C.L. Terng, Critical point theory and submanifold 
geometry, Lecture Notes 

in Math. {\bf 1353}, Springer, Berlin, 1988.

%
%
\noindent
[T] C.L. Terng, 
Proper Fredholm submanifolds of Hilbert space, J. Differential Geometry {\bf 29} 

(1989) 9-47.

\noindent
[TT] C.L. Terng and G. Thorbergsson, 
Submanifold geometry in symmetric spaces, 
J. Differential 

Geometry {\bf 42} (1995) 665-718.

\noindent
[W] B. White, 
The space of $m$-dimensional surfaces that are stationary for a parametric elliptic 

integrand, Indiana Univ. Math. J. {\bf 36} (1987) 567-602.


}


\vspace{0.25truecm}

\rightline{Department of Mathematics, Faculty of Science, }
\rightline{Tokyo University of Science}
\rightline{1-3 Kagurazaka Shinjuku-ku,}
\rightline{Tokyo 162-8601, Japan}
\rightline{(e-mail: koike@ma.kagu.tus.ac.jp)}

\end{document}